\documentclass[11pt]{article}

\usepackage{fullwidth}
\usepackage[top=1.5in, bottom=1.5in, left=1in, right=1in]{geometry}
\usepackage{amsgen}
\usepackage{amsmath}
\usepackage{amstext}
\usepackage{amsbsy}
\usepackage{amsopn}
\usepackage{amsfonts}
\usepackage{amssymb}
\usepackage{eepic}
\usepackage{comment}

\usepackage{tikz-cd}

\usepackage{tikz}
\usetikzlibrary{arrows,automata}
\usetikzlibrary{patterns}
\usetikzlibrary{shapes.geometric,calc}
\tikzstyle{arrow}=[thick, ->, >=stealth]
\tikzset{ edge/.style={->,> = latex'} }
\usetikzlibrary{calc,through,backgrounds,shapes,matrix}
\usepackage{graphicx}
\usepackage{mathrsfs}
\usepackage{faktor}
\usepackage{amsthm}
\numberwithin{equation}{section}

\usepackage[dvipsnames]{xcolor}
\usepackage[colorlinks, linkcolor=BrickRed,citecolor=Violet]{hyperref}
\usepackage[noabbrev,nameinlink]{cleveref}

\def\Box{\square}

\def\tra#1{\smash{\mathop{\mid\kern
-1pt\joinrel\relbar\joinrel\relbar}\limits^{*}_{#1}}}
\def\longtra#1{\smash{\mathop{\mid\kern
-1pt\joinrel\relbar\joinrel\relbar\joinrel\relbar}\limits^{*}_{#1}}}
\def\vlongtra#1{\smash{\mathop{\mid\kern
-1pt\joinrel\relbar\joinrel\relbar\joinrel\relbar\joinrel\relbar}\limits^{*}_{#1}}}
\def\vvlongtra#1{\smash{\mathop{\mid\kern
-1pt\joinrel\relbar\joinrel\relbar\joinrel\relbar\joinrel\relbar\joinrel\relbar}\limits^{*}_{#1}}}
\def\vvvlongtra#1{\smash{\mathop{\mid\kern
-1pt\joinrel\relbar\joinrel\relbar\joinrel\relbar\joinrel\relbar\joinrel\relbar\joinrel\relbar}\limits^{*}_{#1}}}
\def\etra#1{\smash{\mathop{\mid\kern
-1pt\joinrel\relbar\joinrel\relbar}\limits_{#1}}}

\def\gst{{\operatorname{gst}}}

\def\C{{\cal{C}}}

\def\bi{\begin{itemize}}
\def\ei{\end{itemize}}
\def\beq{\begin{equation}}
\def\eeq{\end{equation}}

\def\GM{\operatorname{GM}}
\def\RLM{\operatorname{RLM}}
\def\eval{\operatorname{Eval}}
\theoremstyle{plain}
\newtheorem{T}{Theorem}[section]
\newcommand{\bt}{\begin{T}}
\newcommand{\et}{\end{T}}
\newcommand{\ftd}{$\square$\end{T}}
\newcommand{\Eval}{\operatorname{Eval}}
\newcommand{\States}{\mathsf{States}}

\newtheorem{Proposition}[T]{Proposition}
\newcommand{\bp}{\begin{Proposition}}
\newcommand{\ep}{\end{Proposition}}
\newcommand{\fpd}{$\square$\end{Proposition}}

\newtheorem{Lemma}[T]{Lemma}
\newcommand{\bl}{\begin{Lemma}}
\newcommand{\el}{\end{Lemma}}
\newcommand{\fld}{$\square$\end{Lemma}}

\newtheorem{Corol}[T]{Corollary}
\newcommand{\bc}{\begin{Corol}}
\newcommand{\ec}{\end{Corol}}
\newcommand{\fcd}{$\square$\end{Corol}}

\newtheorem{Result}[T]{Result}
\newcommand{\br}{\begin{Result}}
\newcommand{\er}{\end{Result}}
\newcommand{\frd}{$\square$\end{Result}}

\theoremstyle{definition}
\newtheorem{Example}[T]{Example}
\newcommand{\be}{\begin{Example}}
\newcommand{\ee}{\end{Example}}

\newtheorem{Problem}[T]{Problem}
\newcommand{\bq}{\begin{Problem}}
\newcommand{\eq}{\end{Problem}}

\newtheorem{Remark}[T]{Remark}
\newcommand{\brm}{\begin{Remark}}
\newcommand{\erm}{\end{Remark}}

\newtheorem{Definition}[T]{Definition}
\newcommand{\bd}{\begin{Definition}}
\newcommand{\ed}{\end{Definition}}

\newtheorem{Construction}[T]{Construction}
\newcommand{\bco}{\begin{Construction}}
\newcommand{\eco}{\end{Construction}}

\newtheorem{Hypothesis}[T]{Hypothesis}
\newcommand{\bh}{\begin{Hypothesis}}
\newcommand{\eh}{\end{Hypothesis}}

\newcommand{\opn}{\operatorname}

\normalbaselines
\usepackage[colorinlistoftodos]{todonotes}

\title{Decidability of Krohn-Rhodes Complexity 1 for Finite Semigroups and Automata}
\author{Stuart Margolis, John Rhodes and Anne Schilling}

\date{\today
\footnote{MSC classes: Primary 20M07, 20M10, 20M20, Secondary 54H15}
}

\begin{document}

\maketitle

\begin{abstract}
When decomposing a finite semigroup into a wreath product of groups and aperiodic semigroups, complexity measures the
minimal number of groups that are needed. Determining an algorithm to compute complexity has been an open problem 
for almost 60 years. The main result of this paper proves decidability of Krohn--Rhodes complexity $c = 1$ of finite 
semigroups and automata. This is achieved by showing the lower bounds in work by Henckell, Rhodes and Steinberg from 2012
is sharp using profinite methods and results of McCammond from 1991 and 2001.
\end{abstract}

%\begin{center}
%    {\bf Abstract}
%\end{center}
%
% We prove that it is decidable if a finite semigroup has Krohn--Rhodes complexity 1. This solves a problem that has been open for over 50 years. 

% ----------------------------------------------------------------
\section{Introduction}
% ----------------------------------------------------------------

We prove that it is decidable if a finite semigroup has Krohn--Rhodes complexity 1. This solves a problem that was open for more than 50 years~\cite{KRannals}. In a subsequent paper \cite{complexityn} we use the results of this paper to prove that Krohn--Rhodes complexity is decidable for all $n\geqslant 0$. From hereon in we write ``complexity'' for Krohn--Rhodes complexity.

We review the basic results and tools of complexity theory. We look at some lower and upper bounds to complexity. For more background in complexity theory see~\cite[Chapters 5-9]{Arbib}, \cite[Chapters 4]{qtheory}, and~\cite{TilsonXI, TilsonXII}. See these 
references for all the definitions and results referred to in this paper. All semigroups in this paper are assumed to be finite.

A transformation semigroup (monoid) (abbreviated ts, (tm)) is a pair $(Q,S)$, where $Q$ is a set and $S$ is a subsemigroup (submonoid) of the 
monoid $\mathsf{PT}(Q)$ of all partial functions on $Q$ acting on the right of $Q$. We identify a semigroup $S$ with its right regular representation $(S^{1},S)$, where $S^{1}$ adjoins an identity if $S$ does not already have one.

We begin with the Krohn--Rhodes Theorem originally known as the Prime Decomposition Theorem for Finite Semigroups. Recall that a semigroup 
is {\em aperiodic} if all its subgroups are trivial. A semigroup $S$ divides a semigroup $T$, written $S \prec T$, if $S$ is a quotient of a subsemigroup 
of $T$. There is an analogous definition of division of ts and tm. 

Let $X=(P,T)$ and $Y=(Q,S)$ be ts (tm). Then the wreath product $X \wr Y$ is the ts (tm) $(P \times Q,W)$, where $W$ is the set 
$\{(f,s) \mid s \in S, f\colon\operatorname{Dom}(s)\rightarrow T\}$ where the action is given by $(p,q)(f,s)= (p(qf),qs)$ if $p \in \operatorname{Dom}(qf)$ and undefined otherwise. See the references above for properties of the wreath product of ts.

\bd
A semigroup $S$ is {\em prime} if whenever $S$ divides a wreath product $U \wr T$, then $S$ divides either $U$ or $T$.
\ed

An important semigroup that arises in the theory is the semigroup $\mathsf{RZ}(2)^{1}$ consisting of two right-zeroes and an identity. It is 
faithfully represented on a two-element set by sending the right-zeroes to the two constant functions and the identity to the identity transformation. 
This semigroup is called the flip-flop because of its connection to the theory of finite-state machines. See~\cite{Arbib} for this connection.

\bt[Krohn--Rhodes 1962]\label{KRT}
Every finite semigroup $S$ divides a wreath product of groups and aperiodic semigroups. One can choose the groups to be simple groups that 
divide $S$. We can choose the aperiodic semigroups to be the flip-flop $\mathsf{RZ}(2)^1$.

Furthermore, $S$ is prime if and only if either $S$ is a simple group or $S$ is a subsemigroup of $\mathsf{RZ}(2)^{1}$.
\et

\bc
A semigroup is aperiodic if and only if it divides a wreath product of copies of $\mathsf{RZ}(2)^{1}$.
\ec

The Krohn--Rhodes Theorem was first proved in the joint thesis of Kenneth Krohn and John Rhodes. It first appeared in print in \cite{PDT}. 
Since then it has appeared in a number of books \cite{Arbib, Eilenberg, Lallement, qtheory}. Theorem \ref{KRT} leads to the definition of 
complexity. Complexity was first defined in print in \cite{KRannals}, but had been defined and developed a number of years earlier.

\bd 
The {\em Krohn--Rhodes complexity} of a finite semigroup $S$ is the least number of groups in any decomposition of $S$ as a divisor of wreath
products of groups and aperiodic semigroups.
\ed

The Krohn--Rhodes Theorem guarantees that every finite semigroup $S$ has finite complexity. We write $Sc$ for the complexity of $S$. Clearly, we 
have $Sc=0$ if and only if $S$ is aperiodic and $Sc=1$ if $S$ is a non-trivial group. 

It is known that the full transformation semigroup on $n$ elements $T_n$ has complexity $n-1$ \cite{Arbib, qtheory}. Therefore, if we define $C_n$ to be the set of semigroups 
of complexity at most $n$, then $C_{n}$ is properly contained in $C_{n+1}$ for all $n \geqslant 0$. We will indicate how to prove that $T_{n}c=n-1$ later 
in the paper.

The computability question can be stated as follows. Given a semigroup $S$, say by its multiplication table, the question is whether $Sc$ can be computed. 
The problem is not immediately decidable since one needs to search the infinite set of all wreath product decompositions of $S$ to find and guarantee that the decomposition is minimal. 

The main theorem of this paper proves that the computability question has a positive 
solution if the complexity of $S$ is at most 1. The paper \cite{complexityn} proves that the computability question has a positive solution for all finite semigroups. We state the main theorem.

\bt \label{MainTh} 
{\bf Decidability of Complexity 1 for Finite Semigroups.} Let $S$ be a group mapping semigroup such that its right letter mapping image $\RLM(S)$ has complexity at most 1. Then there is an effectively computable aperiodic transformation semigroup $(Q,A)$ such that the complexity of $S$ is 1 if and only if there is a complete flow $F\colon Q \rightarrow \operatorname{Rh}_{B}(G)$ from $Q$ to the Rhodes lattice $\operatorname{Rh}_{B}(G)$ such that for all $q \in Q$, $qF$ does not contain a contradiction.
\et

There is a large amount of material that needs to be discussed in order to even understand the statement of Theorem \ref{MainTh}. For the convenience of the reader, we summarize all the basic results in finite semigroup theory and further background that is needed for the proof of Theorem \ref{MainTh} in Sections~\ref{section.bounds} to~\ref{section.tools} of this paper. In Section~\ref{section.outline} we give an outline of the proof of Theorem \ref{MainTh}. The rest of the paper is dedicated to the proof of Theorem \ref{MainTh}.

% ----------------------------------------------------------------
\section{Upper and Lower Bounds to Complexity}
\label{section.bounds}
% ----------------------------------------------------------------

% ----------------------------------------------------------------
\subsection{Upper Bounds}

In this section, we discuss computable upper and lower bounds for complexity. Any decomposition of a finite semigroup $S$ as a divisor of a wreath 
product of groups and aperiodic semigroups gives an upper bound to $Sc$. In particular, any proof of the Krohn--Rhodes Theorem gives an upper 
bound to complexity.

The following result of Rhodes and Tilson~\cite{TilsonXI} gives an upper bound based on Green-Rees theory of finite semigroups. We assume familiarity 
with Green--Rees theory.

\bt[The Depth Decomposition Theorem]\label{Depth}
Define the depth $S\delta$ of a finite semigroup to be the longest chain of non-aperiodic $\mathcal{J}$-classes in $S$. Then $Sc \leqslant S\delta$.
\et

The proof of Theorem~\ref{Depth} gives an explicit decomposition of a semigroup $S$, where there are $S\delta$ non-trivial groups. It is well known that 
for the full transformation monoid $T_{n}$, we have $T_{n}\delta=n-1$ and thus Theorem \ref{Depth} implies that $T_{n}c \leqslant n-1$. As mentioned in 
the last section, we actually have equality in this case. On the other hand, if $\operatorname{SIS}_n$ is the symmetric inverse semigroup on $n$ elements consisting of 
all partial bijections on an $n$-element set, then $\operatorname{SIS}_{n}\delta=n-1$ and  $\operatorname{SIS}_{n}c\leqslant 1$ as it is known that any inverse semigroup 
divides $S \wr G$, where $G$ is a group and $S$ is a semilattice. We will describe this decomposition for $\operatorname{SIS}_{n}$ in more detail in further sections.

While Theorem~\ref{Depth} depends on the length of chains of images of elements in a ts, there is another upper bound that is based on sizes of kernel 
classes of elements. Recall that if $s$ is a partial function on a set $Q$, then the kernel of $s$, denoted by $\operatorname{ker}(s)$, is the partition on 
the domain of $s$ given by $q\operatorname{ker}(s) q'$ if and only if $qs=q's$. The {\em degree} of $s$ denoted by $sd$ is the maximal cardinality of a partition 
class of $\operatorname{ker}(s)$. If $X=(Q,S)$ is a ts, then its degree $Xd$ is the maximal degree of any element $s \in S$. For example, 
$(Q,\operatorname{SIS}_{Q})d =1$ for any non-empty set $Q$. See \cite{Tilsonnumber} for a proof of the following theorem. See~\cite{cremona, deg2part2} 
for a detailed study of semigroups of degree 2.

\bt
Let $X=(Q,S)$ be a ts. Then $Sc \leqslant Xd$.
\et

This gives another proof that the symmetric inverse semigroup $\operatorname{SIS}_{Q}$ has complexity at most 1 for any set $Q$. On the other hand, if $Q$ 
is an $n$-element set, and $S$ is the semigroup of constant functions on $Q$, then $(Q,S)d=n$, but $Sc=0$.

We now look at an upper bound defined via chains of morphisms.

\bd\label{morphdef}
\mbox{}
\begin{enumerate}
\item{A morphism $f \colon S \rightarrow T$ is {\it aperiodic} if it is 1-1 on subgroups of $S$.  This is equivalent to having the inverse image of each 
idempotent in $\operatorname{Im}(S)$ being an aperiodic subsemigroup of $S$.}
\item{A morphism $S \rightarrow T$ is $\mathcal{L}'$ if, whenever $x,y$ are regular elements of $S$, then $xf=yf$ implies that 
$x\mathcal{L}y$.}
\end{enumerate}
\ed

The next theorem gives the connection between these morphisms and complexity theory. See also Remark~\ref{rlmrmk} and 
Lemma~\ref{rlmlem} below for their use in so called semi-local theory~\cite[Chapter 7-8]{Arbib} and~\cite[Chapter 4]{qtheory} as a basis for flow theory.

\bt
\mbox{}
\begin{enumerate}
\item{Let $A$ be an aperiodic semigroup and $S$ a semigroup. Then the projection $A \wr S \rightarrow S$ is an aperiodic morphism.}
\item{Let $G$ be a group and $S$ be a semigroup. Then the projection $G \wr S \rightarrow S$ is an $\mathcal{L}'$ morphism.}
\end{enumerate}
\et

We have a partial (much deeper!) converse. See \cite[Chapter 7-8]{Arbib} and \cite{qtheory} for the proofs of the assertions in the following Theorem.

\bt\label{morphcomp}
\mbox{}
\begin{enumerate}
\item{{\bf The Fundamental Lemma of Complexity (Rhodes 1969).} Let $f \colon S \rightarrow T$ be a surjective aperiodic morphism. Then $Sc = Tc$.}
\item{Let $f \colon S \rightarrow T$ be a surjective $\mathcal{L}'$ morphism. Then $Sc \leqslant 1+Tc$.}
\end{enumerate}
\et

\brm
The decidability of complexity can be reduced to the question of whether $Sc=Tc$ or $Sc=1+Tc$ when $f\colon S \rightarrow T$ is a surjective 
$\mathcal{L}'$ morphism.
\erm

\bd
\mbox{}
\begin{enumerate}
\item{Let $S$ be a finite semigroup. An $Ap-\mathcal{L}'$ chain is a factorization of the trivial morphism $\tau_{S} \colon S \rightarrow 1$ of the form
$$\tau_{S}=\alpha_{0}\beta_{1}\cdots \beta_{n}\alpha_{n+1},$$
where each $\alpha_{i}$ is an aperiodic morphism and each $\beta_{j}$ is a non-aperiodic $\mathcal{L}'$ morphism.}
\item{Let $S\theta$ be the minimal $n$ for which $S$ has an $Ap-\mathcal{L}'$ factorization.}
\end{enumerate}
\ed

We note that $S{\theta}$ is computable for a finite semigroup since it is decidable if a congruence on a finite semigroup defines an aperiodic morphism or an $\mathcal{L}'$ morphism.

\bt \label{Morchains}
Let $S$ be a finite semigroup.
\begin{enumerate}
\item{$Sc \leqslant S\theta$.}
\item{If $S$ is completely regular, then $Sc = S\theta$.  Therefore complexity is decidable for completely regular semigroups. See \cite[Chapter 9]{Arbib}.}
\item{$\operatorname{SIS}(n)\theta =n-1$. This shows that there are examples of finite semigroups $S$ such that $S\theta$ is arbitrarily bigger than $Sc$.}
\item{Let $S\widehat{\theta}= \operatorname{Min}\{T\theta \mid S$ {\text is a quotient of} $T\}.$ Then $Sc = S\widehat{\theta}.$}
\end{enumerate}
\et

We note that part 4 of Theorem \ref{Morchains} requires a search though an infinite set and thus does not prove that complexity is computable.

%--------------------------------------------------------------------------------
\subsection{Lower Bounds}

Lower bounds to complexity are more difficult to obtain. If you have a decomposition of a semigroup using $n$ groups with $n>0$, how do 
you prove that there are none that use strictly fewer than $n$ groups? We give one lower bound based on the concept of type I and type II 
subsemigroups and subsemigroup chains. 

\bd
Let $S$ be a finite semigroup.
\begin{enumerate}
\item{A subsemigroup $T$ of $S$ is an {\em absolute type I subsemigroup} if it is generated by a chain of $\mathcal{L}$-classes of $T$.}
\item{The {\em type II subsemigroup $S_{II}$} of $S$ is the smallest subsemigroup of $S$ containing all idempotents and closed under weak conjugation: 
if $xyx = x$ for $x,y \in S$, then $xS_{II}y \cup yS_{II}x \subseteq S_{II}$.}
\end{enumerate}
\ed

\brm
A very important theorem of Ash~\cite{Ash} confirmed a conjecture of Rhodes and shows that if $S$ is a semigroup then $S_{II}$ consists of all 
elements of $S$ that are related to the identity element under all relational morphisms from $S$ to a group. This latter subsemigroup is also 
called the group kernel of $S$. Ash's result showing that the group kernel is precisely the type II subsemigroup of a semigroup $S$ shows that 
membership in the group kernel is decidable. This plays a crucial role in finite semigroup theory and in particular in the proof of decidability of complexity.
\erm

Here is the motivation for type I and type II subsemigroups:

\bt
Let $S$ be a semigroup.
\begin{enumerate}
\item{Assume that $S$ divides $T \wr G \wr A$ for some semigroup $T$, group $G$ and aperiodic semigroup $A$. If S' is 
an absolute type I subsemigroup of $S$, then
$S'$ divides $T \wr A' \wr G'$, where $A'$ is aperiodic and $G'$ is a group.}
\item{Assume that a semigroup $S$ divides $T \wr G$, where $T$ is a semigroup and $G$ is a group. Then $S_{II}$ divides $T$.}
\end{enumerate}
\et

\bc
Let $S$ be a semigroup with complexity $n>0$. Let $T$ be an absolute type I subsemigroup of $S$. Then $T_{II}c < Sc$.
\ec

\bt[Rhodes and Tilson 1972 \cite{lowerbounds2}] 
Let $S$ be a semigroup. Define $Sl$ to be the largest integer $k$ such that there is a chain of subsemigroups
$$ S \geqslant T_{1} > (T_{1})_{II} > T_{2} > (T_{2})_{II} > \cdots >T_{k} >(T_{k})_{II},$$
where each $T_{i}$ is a non-aperiodic absolute type I semigroup. Then $Sl \leqslant Sc$.
\et

It is known that the full transformation semigroup on $n$ elements denoted by $T_{n}$ is generated by its group of units  $\operatorname{Sym}_{n}$ and any element of rank $n-1$ and is thus an absolute type I semigroup. 
Furthermore $T_{n} - \operatorname{Sym}_{n}$ is idempotent generated and it follows that $(T_{n})_{II} = (T_{n} - \operatorname{Sym}_{n}) \cup \{1\}$. Since this semigroup 
contains a copy of $T_{n-1}$ it follows by induction that $T_{n}l= n-1 \leqslant T_{n}c$. This along with the upper bound for $T_{n}c$ that we 
obtained above from the Depth Decomposition Theorem proves that $T_{n}c=n-1$. In the case of the symmetric inverse semigroup 
$\operatorname{SIS}_{n}$, one sees that it is an absolute type I semigroup and that its type II subsemigroup is its semilattice of idempotents. 
It follows that $Sc=Sl \leqslant 1$ for all $n$.

In 1977 \cite{KernelSystems} the first example of a semigroup $S$ for which $Sl < Sc$ was published. That is, $Sl$ is in general a proper 
lower bound to $Sc$. In~\cite{TypeIIfalls}, examples of semigroups $S_{n}$ with $n>0$ such that $S_{n}l=1$ but $S_{n}c=n$ are constructed. 
These examples led to the construction of the lower bound in~\cite{Trans} that turns out to be perfect for complexity 1 as shown in this paper and a
modified version computes arbitrary complexity \cite{complexityn}. On the other hand, $Sl$ is the maximal local complexity function in the sense 
of \cite{localcomplex}. 

The next result shows how to compute the maximal $\mathcal{L}'$ congruence on a semigroup.

\bd
Let $S$ be a semigroup. Let $\mathscr{L}$ be the set of regular $\mathcal{L}$-classes of $S$. Let $l \in \mathscr{L},s \in S$. Then $S$ acts by 
partial functions on $\mathscr{L}$ by $l\cdot s =ls$ if $ls$ is in the same $\mathcal{J}$-class as $l$ and undefined otherwise. Let $S^{\mathcal{L}'}$ 
be the image of $S$ under this action.
\ed

\bt\label{maxLprime}
Let $S$ be a finite semigroup. Then:
\begin{enumerate}
\item{The morphism $\sigma_{L} \colon S \rightarrow S^{\mathcal{L}'}$ is a surjective $\mathcal{L}'$ morphism.}
\item{If $f \colon S \rightarrow T$ is a surjective $\mathcal{L}'$ morphism, then there is a unique surjective $\mathcal{L}'$ morphism 
$g \colon T \rightarrow S^{\mathcal{L}'}$ such that $fg=\sigma_{L}$.}
\item{The assignment of $S$ to $S^{\mathcal{L}'}$ is the object part of a functor on the category of finite semigroups.}
\end{enumerate}
\et

The type I-type II lower bound was used by Tilson and Rhodes~\cite{lowerbounds2}  to show that it is decidable if there exists an aperiodic semigroup 
$A$ and a group $G$ such that  $S$ divides  $A \wr G$. Karnofsky and Rhodes~\cite{complex1/2} used the $Ap-\mathcal{L}'$ upper bound to prove 
that it is decidable if there exists an aperiodic semigroup $A$ and a group $G$ such that  $S$ divides  $G \wr A$. Neither of these results generalize 
to the case of deciding complexity 1, that is if there exist aperiodic semigroups $A_{1},A_{2}$ and a group $G$ such that $S$ divides $A_{1}\wr G \wr A_{2}$.

\bt[Rhodes--Tilson 1972~\cite{{lowerbounds2}}]
Let $S$ be a semigroup. Then $S$ divides $A\wr G$ where $A$ is aperiodic and $G$ is a group if and only if $S_{II}$ 
is aperiodic. It follows that $Sc=Sl$. This condition is decidable.
\et

\bt[Rhodes--Karnofsky 1979~\cite{complex1/2}] Let $S$ be a semigroup. Then $S$ divides $G\wr A$, where $G$ is a group and $A$ is aperiodic if and only if 
$S^{\mathcal{L}'}$ is aperiodic. It follows that $Sc=S\theta$. This condition is decidable.
\et

%--------------------------------------------------------------------------------------------------
\section{Applying the Ideas and Methods of Complexity Theory  to Inverse Semigroups}
\label{section.inverse semigroups}
%--------------------------------------------------------------------------------------------------

Inverse semigroups play an important role in their own right in many parts of mathematics. This is because they have faithful representations 
as ts of partial bijections on sets. 
In his Ph.D. thesis in 1968, Tilson proved that the symmetric inverse semigroup on a set $Q$, $\operatorname{SIS}_{Q}$ divides 
$(\{0,1\},\{0,1\}) \wr (Q,\operatorname{Sym}_{Q})$, 
the wreath product of the symmetric group on $Q$ and the two-element semilattice $\{0,1\}$. By the Preston--Wagner Theorem~\cite{Lawson},
every inverse semigroup $S$ is isomorphic to a subsemigroup of $\operatorname{SIS}_{S}$ and it follows that every inverse semigroup has complexity at most 1.

In this section, we show that many of the important concepts of inverse semigroup theory including $E$-unitary inverse semigroups and fundamental
 inverse semigroups have natural interpretations in terms of complexity theory. Congruences on ts whose quotients act by partial bijections play a crucial 
 role in complexity theory. It is useful to explicitly explain the complexity theoretic aspects of inverse semigroup theory.

\bd
We give a complexity theoretic definition of the important class of $E$-unitary semigroups. An inverse semigroup $S$ is {\em $E$-unitary} if the 
morphism $\sigma_{Gp} \colon S \rightarrow G$ from $S$ to its maximal group image is an aperiodic morphism.  This is equivalent to this morphism 
being idempotent-pure, that is, the idempotents form a congruence class of this morphism. This is one of the standard definitions of $E$-unitary semigroups.
\ed

\brm
It follows from this definition that if $S$ is an $E$-unitary inverse semigroup, then $\tau_{S}=\sigma_{Gp}\tau_{G}$ is an $Ap-\mathcal{L}'$ 
factorization of the trivial morphism on $S$. Consequently, $S\theta \leqslant 1$. 
\erm

\bt[O'Carroll~\cite{OCarroll}]
\label{OCaroll} 
Let $S$ be an inverse semigroup. Then there is a semilattice $T$ such that $S$ is a subsemigroup of $T \wr G$, where $G$ is the maximal group 
image of $S$ if and only if $S$ is $E$-unitary.
\et

\bt[McAlister \cite{McAlisterPstuff}]
 \label{Ecov}
Let $S$ be an inverse semigroup. Then there is an $E$-unitary inverse semigroup and a surjective idempotent-separating morphism 
$T \rightarrow S$. Consequently $S\widehat{\theta}= Sc \leqslant 1$.
\et

Any inverse semigroup with 0, (for example, $\mathsf{SIS_{Q}}$) has maximal group image the trivial group. 
Thus an inverse semigroup with 0 is $E$-unitary and thus embeds into the semidirect product of a semilattice and a group by Theorem \ref{OCaroll} 
if and only if it is a semilattice. Therefore the result that every inverse semigroup divides a semidirect product of a semilattice and a group can not 
be replaced by an embedding theorem into such a semidirect product. 

Let $T$ be an $E$-unitary cover of $\operatorname{SIS}_{Q}$. Its maximal group image (remember all semigroups in this paper are finite) is isomorphic to
its minimal ideal $G$. The group $G$ maps onto the 0 of the $\operatorname{SIS}_{Q}$ by the covering morphism. We interpret the group $G$ as resolving 
the ``singularity" of the 0 of $\operatorname{SIS}_{Q}$ that is preventing the embedding of $\operatorname{SIS}_{Q}$ as a subsemigroup of a semidirect product 
of a semilattice and a group by Theorem~\ref{OCaroll}. Expanding semigroups by removing such singularities is an important tool in semigroup 
theory and in the proof of computability of complexity of finite semigroups.

We give another approach to the complexity theory of inverse semigroups via a direct product decomposition. This will generalize via the 
Presentation Lemma and the Theory of Flows to arbitrary semigroups.

Another important class of inverse semigroups is the class of fundamental inverse semigroups. An inverse semigroup is {\em fundamental} if the largest 
congruence contained in $\mathcal{H}$ is the trivial congruence. There is a maximal congruence contained in $\mathcal{H}$ on any inverse 
semigroup $S$. Its quotient is the maximal fundamental image of $S$. 

From Theorem \ref{maxLprime}  the maximal $\mathcal{L}'$ congruence on an inverse semigroup is given by the kernel of the morphism from 
$S$ to $S^{\mathcal{L}'}$. Clearly, on a regular semigroup this is the maximal congruence contained in $\mathcal{L}$. In an inverse semigroup 
we also have that this is the maximal congruence contained in $\mathcal{R}$ by inversion. Thus this is the largest congruence contained in 
$\mathcal{H}$ for an inverse semigroup. The following theorem follows immediately. This gives a complexity theoretic interpretation of 
fundamental inverse semigroups. The reader familiar with inverse semigroup can easily check that the representation defining $S^{\mathcal{L}'}$ 
is equivalent to the Munn representation~\cite{Lawson} that is the usual way to compute the fundamental image of an inverse semigroup.

\bt
Let $S$ be an inverse semigroup. Then $S^{\mathcal{L}'}$ is the maximal fundamental image of $S$. Consequently $S$ is a fundamental inverse 
semigroup if and only if $S$ is isomorphic to $S^{\mathcal{L}'}$.
\et

The next important theorem of McAlister and Reilly \cite{McAlReilly} shows that an inverse semigroup divides the direct product of a group and 
its maximal fundamental image. This is a special case of the Presentation Lemma Decomposition Theorem discussed below.

\bt[McAlister--Reilly 1975~\cite{McAlReilly}]
\label{McRe}
Let $S$ be an inverse semigroup. Then $S$ divides \linebreak
$(H \wr \operatorname{Sym}_{\mathscr{L}}) \times S^{\mathcal{L}'}$, where $H$ is the direct 
product of the maximal subgroups of $S$, one for each $\mathcal{D}$-class, and $\mathscr{L}$ is the set of $\mathcal{L}$-classes of $S$.
\et

%-----------------------------------------------------------------------------------------------------
\section{Tools Used for the Proof of Computability of Complexity}
\label{section.tools}
%-----------------------------------------------------------------------------------------------------

%-----------------------------------------------------------------------------------------------------
\subsection{Transitive Representations of Finite Semigroups}

We look at the tools and ideas used in the proof of the main theorem by reducing the problem to semigroups faithfully represented by partial 
functions acting transitively on a set.  We first review this theory.

\bd
A semigroup $S$ is right (left, bi-) transitive if $S$ has a faithful transitive right (left, right and left) action by partial functions on some set $X$.
\ed

\brm
In the literature these are called Right Mapping, Left Mapping and Generalized Group Mapping semigroups respectively~\cite{Arbib, qtheory}.
\erm

The most important example of a transitive action is the Sch\"utzenberger representation $\operatorname{Sch}(R)$ on an $\mathcal{R}$-class 
$R$ of $S$.

For $r \in R, s \in S$, define 
\[
r\cdot s =	\begin{cases}
      rs & \text{if}\ rs \in R ,\\
      \text{undefined} & \text{otherwise.}
    \end{cases}
\]
There is also the dual notion of the left Sch\"utzenberger representation $\operatorname{Sch}(L)$ on an $\mathcal{L}$-class $L$ of $S$.

\begin{Lemma}
The semigroup $S$ is right (left, bi-) transitive if and only if $S$ has a unique $0$-minimal regular ideal $I(S) \approx M^{0}(G,A,B,C)$ such that 
the right (left, right and left) Sch\"utzenberger representation of $S$ is faithful on any $\mathcal{R}, (\mathcal{L}, \mathcal{R}$ and 
$\mathcal{L})$-class $R (L, R$ and $L)$ in $I(S) -\{0\}$.  
\end{Lemma}

\bd
A bi-transitive semigroup is called {\em Group Mapping} written $\operatorname{GM}$ if the group $G$ in $I(S)$ is non-trivial.
\ed

\brm
\label{rlmrmk}
We make the following remarks. 

\begin{enumerate}
\item{We can identify a fixed $\mathcal{R}$-class in $I(S)$ with the set $G \times B$. We obtain a faithful transformation semigroup 
$(G \times B,S)$.}
\item{The action of $S$ on $G \times B$ induces an action of $S$ on $B$.  The faithful image of this action is called the 
{\em Right Letter Mapping} (written $\RLM(S)$) image of $S$.  We thus obtain a faithful transformation semigroup $(B,\RLM(S))$.}
\end{enumerate}
\erm

See Chapters 7-8 of \cite{Arbib}, or Chapter 4 of \cite{qtheory} for proofs of the next Lemma. 

\begin{Lemma}\label{rlmlem}
Let $S$ be a $\GM$ semigroup.
\begin{enumerate}
\item{$(G \times B,S)$ embeds into $G \wr (B,\RLM(S))$. Therefore, $Sc \leqslant 1+\RLM(S)c$}
\item{$S^{\mathcal{L}'} = \RLM(S)$.}
\item{For every semigroup $T$, there is a $\GM$ quotient semigroup $S$ such that $Tc=Sc$.} 
\end{enumerate}
\end{Lemma}

The following important theorem explains why flow theory is based upon $\GM$ and $\RLM$ semigroups. It follows immediately from Lemma \ref{rlmlem} and induction on the cardinality of $S$, since 
$|\RLM(S)| < |S|$ when $S$ is a $\GM$ semigroup.

\bt
The question of decidability of complexity can be reduced to the case that $S$ is a $\GM$ semigroup and whether $Sc=\RLM(S)c$ or 
$Sc = 1+\RLM(S)c$.
\et

%---------------------------------------------------------------------------------------
\subsection{The Set-Partition Lattice and the Rhodes Lattice of a Group-Mapping Semigroup }
\label{PLFlows}

The Presentation Lemma gives a necessary and sufficient condition for a $\GM$ semigroup $S$ to have the property that $Sc = \RLM(S)c$. 
The idea goes back to Tilson's proof that complexity is computable for semigroups with at most 2 non-zero $\mathcal{J}$-classes \cite{2J}. 
There are various versions in the literature. See~\cite{AHNR.1995}, \cite[Section 4.14]{qtheory} for example. 
In~\cite{Trans} the notion of flows was defined, developed and proved to give an alternative way of describing the Presentation Lemma. 
We will use flows in this paper. See~\cite{FlowsI, FlowsII} for the theory of flows over aperiodic semigroups. 

Let $S$ be a $\operatorname{GM}$ semigroup with $0$-minimal ideal $M^{0}(G,A,B,C)$. The definition of flow gives a map from the state set 
of an automaton to a lattice associated with $S$. In the literature, there are two such lattices. We define these lattices and give the connections 
between them. The first is the set-partition lattice $\operatorname{SP}(G \times B)$. This is the lattice 
whose elements are all pairs $(Y,\Pi)$, where $Y$ is a subset of $G\times B$ and $\Pi$ is a partition on $Y$. Here $(Y,\Pi) \leqslant (Z,\Theta)$ if 
$Y \subseteq Z$ and for all $y \in Y$, the $\Pi$ class of $y$ is contained in the $\Theta$ class of $y$. 

The second lattice is the Rhodes lattice $\operatorname{Rh}_{B}(G)$. We review the basics. For more details see~\cite{AmigoDowling}. 
Let $G$ be a finite group and $B$ a finite set. A partial partition on $B$ is a partition $\Pi$ on a subset $I$ of $B$. We also consider the 
collection of all functions $F(B,G)$, $f \colon I \rightarrow G$ from subsets $I$ of $B$ to $G$. The group $G$ acts on the left of $F(B,G)$ by $(gf)(b) = gf(b)$
for $f \in F(B,G), g \in G, b \in \operatorname{Dom}(f)$. An element $\{gf \mid f \colon I \rightarrow G, I \subseteq X, g \in G\}$ of the quotient set 
$F(B,G)/G$ is called a cross-section 
with domain $I$. It should be thought of as a projectification of a cross-section of the projection from $I$ to $B$ in the usual topological sense. An
SPC (Subset, Partition, Cross-section) over $G$ is a triple $(I,\Pi,f)$, where $I$ is a subset of $B$, $\Pi$ is a partition of $I$ and $f$ is a collection of 
cross-sections one for each $\Pi$-class $\pi$ with domain $\pi$. If the classes of $\Pi$ are $\{\pi_{1},\pi_{2}, \ldots, \pi_{k}\}$, then we sometimes write 
$\{(\pi_{1},f_{1}), \ldots, (\pi_{k},f_{k})\}$, where $f_{i}$ is the cross-section associated to $\pi_{i}$. For brevity we denote this set of cross-sections by 
$[f]_{\pi}$.  We let $\operatorname{Rh}_{B}(G)$ denote the set of all SPCs on $B$ over the group $G$ union a new element 
$\Longrightarrow\Longleftarrow$  that we call {\em contradiction} and is the top element of the lattice structure on $\operatorname{Rh}_{B}(G)$.
Contradiction occurs because the join of two SPCs need not exist. In this case we say that the contradiction is their join.

The partial order on $\operatorname{Rh}_{B}(G)$ is defined as follows. We have $(I,\pi,[f]_{\pi}) \leqslant (J,\tau,[h]_{\tau})$ if:
\begin{enumerate}
\item
{$I \subseteq J$;} 
\item
{Every block of $\pi$ is contained in a (necessarily unique) block of $\tau$;} 
\item{if the $\pi$-class $\pi_{i}$ is a subset of the $\tau$-class $\tau_{j}$, then the restriction of $h$ to $\pi_{i}$ equals $f$ restricted to $\pi_{i}$ as elements of $F(B,G)/G$. That is, $[h|_{\pi_{i}}] = [f|_{\pi_{i}}]$.}
\end{enumerate}

See \cite[Section 3]{AmigoDowling} for the definition of the lattice structure on $\operatorname{Rh}_{B}(G)$. The underlying set of the Rhodes lattice 
$\operatorname{Rh}_{B}(G)$ minus the contradiction is the set underlying the Dowling lattice on the same set and group. The Dowling lattice 
has a different partial order. For the connection between Rhodes lattices and Dowling lattices see \cite{AmigoDowling}.

To be consistent with the paper \cite{Trans} we will use flows with respect to the set-partition lattice $\operatorname{SP}(G \times B)$. Since Rhodes lattices and flows to them also appear in a number of places in the literature, we give the connections between the lattices referring to the literature for complete proofs.

We note that $\operatorname{SP}(G \times B)$ is isomorphic to the Rhodes lattice $\operatorname{Rh}_{G\times B}(1)$ of the trivial group over the 
set $G \times B$. We need only note that a cross-section to the trivial group is a partial constant function to the identity and can be omitted, leaving 
us with a set-partition pair. There are no contradictions for Rhodes lattices over the trivial group, and in this case the top element is the pair 
$(G \times B,(G \times B)^{2})$. Despite this, we prefer to use the notation $\operatorname{SP}(G \times B)$ instead of $\operatorname{Rh}_{G\times B}(1)$.

Conversely, we can find a copy of the meet-semilattice of $\operatorname{Rh}_{B}(G)$ as a meet subsemilattice of $\operatorname{SP}(G\times B)$. 
We begin with the following important definition.

\bd
A subset $X$ of $G \times B$ is a cross-section if whenever $(g,b),(h,b) \in X$, then $g = h$. That is, $X$ defines a cross-section of the projection 
$\theta \colon X \rightarrow B$. Equivalently, $X^{\rho} \subseteq B \times G$, the reverse of $X$, is the graph of a partial function 
$f_{X} \colon B \rightarrow G$.  An  element $(Y,\Pi) \in \operatorname{SP}(G \times B)$ is a cross-section if  every partition class $\pi$ of $\Pi$ 
is a cross-section.  
\ed

From the semigroup point of view, a cross-section is a partial transversal of the $\mathcal{H}$-classes of the distinguished $\mathcal{R}$-class, 
$R = G \times B$ of a $\GM$ semigroup. That is, the $\mathcal{H}$-classes of $R$ are indexed by $B$ and a cross-section picks at most one 
element from each $\mathcal{H}$-class. 

An element $(Y,\Pi)$ that is not a cross-section is called a {\em contradiction}. That is, $(Y,\Pi)$ is a contradiction if some $\Pi$-class $\pi$ contains 
two elements $(g,b),(h,b)$ with $g \neq h$. We note that the set of cross-sections is a meet subsemilattice of $\operatorname{SP}(G \times B)$ and 
the set of contradictions is a join subsemilattice of $\operatorname{SP}(G \times B)$.

We will identify $B$ with the subset $\{(1,b)\mid b \in B\}$ of $G \times B$, which as above we think of as a system of representatives of the 
$\mathcal{H}$-classes of the distinguished $\mathcal{R}$-class. Note then that $G\times B$ is the free left $G$-act on 
$B$ under the action $g(h,b)=(gh,b)$. This action extends to subsets and partitions of $G\times B$ . Thus $\operatorname{SP}(G\times B)$ 
is a left $G$-act. An element $(Y,\Pi)$ is {\em invariant} if $G(Y,\Pi) = (Y,\Pi)$. It is easy to see that $(Y,\Pi)$ is invariant if and only if:

\begin{enumerate}
 \item{$Y=G \times B'$ for some subset $B'$ of $B$.}  

 \item{For each $\Pi$-class $\pi$, $G\pi \subseteq \Pi$ and is a partition of $G\times B''$ where $B''\subseteq B'$.}
\end{enumerate}

Thus $(Y,\Pi)$ is invariant if and only if $Y=G\times B'$ for some subset $B'$ of $B$ and there is a partition $B_{1},\ldots, B_{n}$ of $B'$ such that for all 
$\Pi$ classes $\pi$, $G\pi$ is a partition of $G\times B_{i}$ for a unique $1 \leqslant i \leqslant n$. 

Let $\operatorname{CS}(G\times B)$ be the set of invariant cross-sections $(Y,\Pi)$ in $\operatorname{SP}(G \times B)$. $\operatorname{CS}(G\times B)$ is a 
meet-subsemilattice of $\operatorname{SP}(G\times B)$. We add a top-element $\Rightarrow\Leftarrow$ to $\operatorname{CS}(G \times B)$. The resulting structure $\operatorname{CS}^{\Rightarrow\Leftarrow}(G \times B)$ is a lattice if we define join to be its join in $\operatorname{CS}(G \times B)$ if it exists and $\Rightarrow\Leftarrow$ otherwise.

The following Proposition is proved in \cite{FlowsI}. It gives the connection between set-partition lattices and Rhodes lattices.

\bp\label{Updown}
The set-partition lattice $\operatorname{SP}(G\times B)$ is isomorphic to the Rhodes lattice $\operatorname{Rh}_{G\times B}(1)$. 
The Rhodes lattice $\operatorname{Rh}_{B}(G)$ is isomorphic to the lattice $\operatorname{CS}^{\Rightarrow\Leftarrow}(G\times B)$. 
\ep

\subsection{Congruences on a Transformation Semigroup and the Tilson Congruence}

A {\em congruence} on a transformation semigroup $(Q,S)$ is an equivalence relation $\approx$ on $Q$ such that if $q\approx q'$ and for $s \in S$, 
and both $qs$ and $q's$ are defined, then $qs\approx q's$. Every $s \in S$ defines a partial function on $\faktor{Q}{\approx}$ by 
$[q]_{\approx}s=[q's]_{\approx}$ if  $q's$ is defined for some $q' \in [q]_{\approx}$. The quotient $\faktor{(Q,S)}{\approx}$ has states 
$\faktor{Q}{\approx}$ and semigroup $T$ which is generated by the action of all $s \in S$ on $\faktor{Q}{\approx}$. We remark that $T$ is 
not necessarily a quotient semigroup of $S$, but is in the case that $(Q,S)$ is a transformation semigroup of total functions. 

A congruence $\approx$ is called {\em injective} if every $s \in S$ defines a partial 1-1 function on $\faktor{Q}{\approx}$. It is easy to see that the 
intersection of injective congruences is injective. Therefore, there is a unique minimal injective congruence $\tau$ on any transformation semigroup 
$(Q,S)$. We call $\tau$ the Tilson congruence on a transformation semigroup because of the following proved in~\cite{Redux}. It is central to the 
theory of flows.
 
\bt\label{redux}
Let $(G \times B,S)$ be a $\GM$ transformation semigroup. Then the minimal injective congruence $\tau$ on $S$ is defined as follows: 
$(g,b) \tau (g',b')$ if and only if there are elements $s,t \in S_{II}$ such that $(g,b)s=(g',b')$ and $(g',b')t=(g,b)$.
\et
 
\brm
\mbox{}
\begin{enumerate}
\item{We can state this theorem by $(g,b)S_{II}=(g',b')S_{II}$. That is, $(g,b)$ and $(g',b')$ define the same ``right coset'' of $S_{II}$ on $G \times B$.}
\item{Theorem \ref{redux} can be used to greatly simplify the proof in \cite{lowerbounds2} for decidability of membership in $S_{II}$ for regular elements 
of an arbitrary semigroup.}
\end{enumerate}
\erm

%---------------------------------------------------------------------------------------
\subsection{Definition of Flows and Their Properties}
Let $(G \times B,S)$ be the transformation semigroup associated to a $\GM$ semigroup $S$ and let $X$ be a generating set for $S$. 
By a deterministic automaton we mean an automaton such that each letter defines a partial function on the state set.

\bd
Let $\mathcal{A}$ be a deterministic automaton with state set $Q$ and alphabet $X$. 
A {\em flow} to the lattice $\operatorname{SP}(G\times B)$ on $\mathcal{A}$ is a function $f \colon Q \rightarrow \operatorname{SP}(G \times B)$ 
such that for each $q \in Q, x \in X$, with $qf =(Y,\pi)$ and $(qx)f = (Z,\theta)$ we have:
\begin{enumerate}
\item{For all $(g,b) \in G \times B$, there is a $q \in Q$ such that $(g,b) \in qf$.}

\item{$Yx \subseteq Z$.}

\item{Multiplication by $x$ considered as an element of $S$ induces a partial 1-1 map from $Y/\pi$ to $Z/\theta$.

That is, for all $y,y' \in Y$ we have $y,y'$ are in a $\pi$-class if and only if $yx,y'x$ are in a $\theta$-class whenever $yx,y'x$ are both defined.}

\item{{\bf The Cross Section Condition:} For all $q \in Q$, $qf$ is a cross-section. That is, for all $g,h \in G, b \in B$ if $(g,b),(h,b) \in qf$ it follows that $g=h$.}

\end{enumerate}
\ed

We use Proposition \ref{Updown} to give a definition of a flow to the Rhodes lattice $\operatorname{Rh}_{B}(G)$.

\bd
Let $\mathcal{A}$ be a deterministic automaton with state set $Q$ and alphabet $X$. 
A {\em flow} to the lattice $\operatorname{Rh}_{B}(G)$ is a flow to $\operatorname{SP}(G\times B)$ such that for each state $q$, 
$qf \in \operatorname{CS}(G\times B)$. That is, $qf$ is a $G$-invariant cross-section.
\ed

\brm
It follows from the original statement of the Presentation Lemma \cite{AHNR.1995} that if $S$ has a flow with respect to some automaton over 
$\operatorname{SP}(G\times B)$, then it has a flow from the same automaton over $\operatorname{Rh}_{B}(G)$. The proof of the equivalence 
of the Presentation Lemma and Flows in~\cite[Section 3]{Trans} works as well for the Presentation Lemma in the sense of \cite{AHNR.1995}. 
\erm

We will use all the terminology and concepts from \cite{Trans}. If $\mathcal{A}$ is an automaton with alphabet $X$ and state set $Q$, recall that its 
completion is the automaton $\mathcal{A}^{\square}$ that adds a sink state $\square$ to $Q$ and declares that $qx = \square$ if $qx$ is not 
defined in $\mathcal{A}$. A flow on $\mathcal{A}$ is a complete flow if on $\mathcal{A}^{\square}$ extending $f$ by letting 
$\square f = (\emptyset, \emptyset)$ the bottom of both the set-partition and Rhodes lattices remains a flow and furthermore, for all 
$(g,b) \in G \times B$, there is a $q\in Q$ such that $(g,b) \in Y$, where $qf=(Y,\Pi)$. All flows in this paper will be complete flows.

Flows are related to the Presentation Lemma \cite{AHNR.1995}, \cite[Section 4.14]{qtheory}. They give a necessary and sufficient condition for 
$Sc =\operatorname{RLM}(S)c$ where $S$ is a $\operatorname{GM}$ semigroup.  See~\cite[Section 3]{Trans} for a proof of the following Theorem.

\bt[The Presentation Lemma-Flow Version]
\label{PLflow}
Let $(G \times B,S)$ be a $GM$ transformation semigroup with $S$ generated by $X$. Let $k> 0$ and assume that $\operatorname{RLM}(S)c = k$.  
Then $Sc = k$ if and only if there is an $X$ automaton $\mathcal{A}$ whose transformation semigroup $T$ has complexity strictly less than $k$ and a 
complete flow $f \colon Q \rightarrow \operatorname{Rh}_{B}(G)$.
\et

 Flows were preceded by the Presentation Lemma \cite{AHNR.1995}, \cite[Section 4.14]{qtheory}. The Presentation Lemma was shown to be equivalent to the existence of an appropriate flow in \cite[Section 3]{Trans}. There are three versions of the Presentation Lemma and its relation to Flow Theory in the literature \cite{AHNR.1995}, \cite[Section 4.14]{qtheory}, \cite{Trans}. These have very different formalizations and it is not clear how to pass from one version to another. The terminology is different as well. For example, the definition of cross-section in each of these references, as well as the one we use in this paper is different.  The following Theorem gives a unified approach to these topics. It summarizes known results in the literature and is meant to emphasize the strong connections between slices, flows and the presentation lemma and the corresponding direct product decomposition. See~\cite{FlowsI} for a proof. We use some of the results that we have proved in previous sections. For background on the derived transformation semigroup $D(\Phi)$ of a relational morphism $\Phi$ and the derived semigroup theorem see \cite{Eilenberg, qtheory}.

\bt \label{uniform}
Let $(G \times B,S)$ be $\GM$ and assume that $\RLM(S)c \leqslant n$. Then the following are equivalent:

\begin{enumerate}

\item{$Sc \leqslant n$.}

\item{There is an aperiodic relational morphism $\Theta:S \rightarrow H\wr T$, where $H$ is a group and $Tc \leqslant n-1$.}

\item{There is a relational morphism $\phi\colon S \rightarrow T$ where $Tc \leqslant n-1$ and such that the derived transformation semigroup $D(\Phi)$ 
is in $Ap*Gp$ where $Ap*Gp$ is the pseudovariety generated by all semidirect products $T \rtimes G$ where $T$ is an aperiodic semigroup and $G$ is a group.}

\item{There is a relational morphism $\phi\colon S \rightarrow T$ where $Tc \leqslant n-1$ such that the Tilson congruence $\tau$ on the derived transformation semigroup $D(\Phi)$ is a cross-section.}

\item{$(G \times B,S)$ admits a flow from a transformation semigroup $(Q,T)$ with $(Q,T)c \leqslant n-1$.}

\item{$S \prec (G \wr (\operatorname{Sym}(B)) \wr T) \times \RLM(S)$ for some transformation semigroup $T$ with $Tc \leqslant n-1$.}

\end{enumerate}

\et

We emphasize the case for complexity 1 and aperiodic flows. 

\bt \label{Apuniform}
Let $(G \times B,S)$ be $\GM$ and assume that $\RLM(S)c \leqslant 1$. Then the following are equivalent:
\begin{enumerate}

\item{$Sc=1$.}

\item{There is an aperiodic relational morphism $\Theta:S \rightarrow H\wr T$, where $H$ is a group and $T$ is aperiodic.}

\item{There is a relational morphism $\phi\colon S \rightarrow T$ where $T$ is aperiodic such that the Derived transformation semigroup $D(\Phi)$ is in $Ap*Gp$.}

\item{There is a relational morphism $\phi\colon S \rightarrow T$ where $T$ is aperiodic such that the Tilson congruence $\tau$ on the Derived transformation semigroup $D(\Phi)$ is a cross-section.}

\item{$(G \times B,S)$ admits an aperiodic flow.}

\item{$S \prec (G \wr (\operatorname{Sym}(B)) \wr T) \times \RLM(S)$ for some aperiodic semigroup $T$.}

\end{enumerate}

\et

We emphasize the direct product decomposition in Theorem~\ref{uniform}.
\bc[The Presentation Lemma Decomposition Theorem]
\label{PLDT}
Let $(G \times B,S)$ be $\GM$ and let $k > 0$.  Then $Sc \leqslant k$ if and only if $\RLM(S)c \leqslant k$ and there is a transformation semigroup 
$T$ with $Tc \leqslant k-1$ such that 
$$ (G \times B,S) \prec (G \wr \operatorname{Sym}_{B} \wr T) \times (B,\RLM(S)).$$
\ec

Note that if $S$ is a $\GM$ inverse semigroup then Theorem \ref{McRe} follows from Corollary \ref{PLDT} if we take $T$ to be the trivial semigroup.  It is immediate from the definition in fact that a $\GM$ inverse semigroup has a flow over the trivial semigroup. See \cite{FlowsI} for a description of semigroups that have flows over the trivial semigroup.

The best way to learn Flow Theory is to study, work through and create examples. A Master List of Examples is being maintained on the ArXiv~\cite{MasterList} and on the Research Gate accounts of the first and second authors. The list is dynamic and will be updated as necessary. We give an important example from that paper called (for historical reasons) $\operatorname{BIRIP}$.

\be

Consider the matrices
\[
	M_2 = \begin{pmatrix}
	1&1\\
	1&0\\
	0&1
	\end{pmatrix}
	\qquad \text{and} \qquad
	M_4 = \begin{pmatrix}
	1&1&0&0\\
	0&1&1&0\\
	0&0&1&1\\
	1&0&0&1\\
	1&0&0&0\\
	0&1&0&0\\
	0&0&1&0\\
	0&0&0&1
	\end{pmatrix}
\]
with indexing sets $\{1',3'\}$ and $\{1,2,3,4\}$, respectively. In addition, consider the partial monomial maps
\[
	x = \begin{cases}
	1' \mapsto 3',\\
	3' \mapsto 1',
	\end{cases}
	\quad
	a=\begin{cases}
	1 \mapsto 4,\\
	2 \mapsto 3,	
	\end{cases}
	\quad
	b = \begin{cases}
	1 \mapsto 3,\\
	3 \mapsto 1,\\
	2 \mapsto 4,\\
	4 \mapsto 2,
	\end{cases}
	\quad \text{and} \quad
	t = \begin{cases}
	1' \mapsto 1,\\
	3' \mapsto -3.
	\end{cases}
\]
It is routine to check that the maps $x,a,b,t$ satisfy the linked equations~\cite[Chapter 7]{Arbib.1968} and~\cite[Chapter 4]{qtheory} with respect 
to the Rees matrix semigroup $S$  with structure matrix, the sum of $M_4$ and $M_2$ over the group $Z_2=\{1,-1\}$.
This implies that the {\em semigroup $\operatorname{BIRIP}$} generated by
$x,a,b,t$ and the Rees matrix semigroup $S$ is a group-mapping semigroup with 
distinguished ideal equal to $S$.

Note that $\mathsf{RLM}(\operatorname{BIRIP})$ has the distinguished ideal $I$ over the trivial group with structure matrix $M_4+M_2$. Since $x,a,b,t$ 
are partial one-to-one maps, the quotient $\mathsf{RLM}(\operatorname{BIRIP})/I$ embeds into a symmetric inverse semigroup. By the Fundamental 
Lemma of Complexity~\ref{morphcomp}, it follows that $\mathsf{RLM}(\operatorname{BIRIP})c=1$.

Consider the aperiodic automaton $\mathcal{A}$ with three states
\begin{center}
\begin{tikzpicture}[auto]
\node (A) at (0, 0) {$\mathbf{1}$};
\node (B) at (4,0) {$\mathbf{2}$};
\node (C) at (8,0) {$\mathbf{3}$};

\path (A) edge [thick,loop left] node {$x$} (A);
\draw[edge,thick] (A) -- (B) node[midway, above] {$t$};
\path (B) edge [thick,loop above] node {$b$} (B);
\draw[edge,thick] (B) -- (C) node[midway, above] {$a$};
\path (C) edge [thick,loop right] node {$a,b$} (C);
\end{tikzpicture}
\end{center}
where in addition $(\alpha,g,\beta) \in I$ acts as the constant map on $\mathbf{1},\mathbf{2},\mathbf{3}$ to
$1$ if $\beta \in \{1',3\}$ and to $3$ if $\beta \in \{1,2,3,4\}$.

We claim that
\[
	\mathbf{1} \times i = \{1',3'\}/\langle 11\rangle, \qquad
	\mathbf{2} \times i = \{1,3\}/\langle 1-1\rangle, \qquad
	\mathbf{3} \times i = \{1,2,3,4\}/\langle 1111\rangle
\]
is a flow (in downstairs notation). Here $\mathbf{1}\times i$ lies in the set of $\mathsf{States}$ , because 
\[
	\mathbf{1}\times i = 1' x^{\omega+\star}.
\]
Similarly, $\mathbf{2}\times i$ lies in $\mathsf{States}$, because
\[
	\mathbf{2} \times i = (\mathbf{1} \times i)t.
\]
Finally, $\mathbf{3} \times i$ lies in $\mathsf{States}$, because 
\[
	\mathbf{3} \times i = 1't(a b^{\omega+\star})^{\omega+\star}.
\]
The reader should check that $\times i$ is indeed a flow.

Hence, by the Presentation Lemma we have $(S,X)c=1$.

\ee

%--------------------------------------------------------------------------------------------------------
\subsection{More tools used in the proof of computability of complexity 1}
\label{section.more tools}

The proof of computability of complexity 1 computes a finite state automaton $\mathcal{A}$ with a state set $Q$ and aperiodic transformation 
semigroup $T$ such that there is either a flow $f \colon Q \rightarrow \operatorname{SP}(G\times B)$ defined on $\mathcal{A}$ or no aperiodic flow exists.
We give a list of the tools, steps and ideas used in the proof.  

\begin{itemize}

\item{{\color{blue} Henckell's Theorem on Pointlike Aperiodic Sets} \cite{Henckell, BenHenckell}  

A subset $Y$ of a semigroup $S$ is {\em pointlike for aperiodic semigroups} if for all relational morphisms $\theta \colon S \rightarrow T$ with $T$ aperiodic, 
there is a $t \in T$ such that $Y \subseteq t\theta^{-1}$.

Henckell proved \cite{Henckell} that the pointlike subsets of $S$ are the smallest subsemigroup of the semigroup of subsets $P(S)$ containing the 
singletons, closed under subset and closed under the operation that sends a set $Z$ to $Z^{\omega}(Z \cup Z^{2} \cup \cdots \cup Z^{\omega})$,
where $Z^{\omega}$ is the unique idempotent in the subsemigroup generated by $Z$. 
\smallskip
It follows that it is decidable if a subset of $S$ is pointlike.}

\end{itemize}

We record more properties of pointlike sets that we will use in this paper. All of these properties hold more generally when studying pointlike sets relative to an arbitrary pseudovariety, but in this paper we only need the case of the pseudovariety of aperiodic semigroups. See \cite{qtheory} for the proofs of the properties noted here.

Note that $\mathsf{Pl_{Ap}}(S)$ satisfies all of the following (not necessarily independent) properties:
\begin{enumerate}
\item (points) $\{s\} \in \mathsf{Pl_{Ap}}(S)$;
\item (subsets) If $Y \subseteq X$ for $X \in \mathsf{Pl_{Ap}}(S)$, then $Y \in \mathsf{Pl_{Ap}}(S)$;
\item (products) If $X,Y \in \mathsf{Pl_{Ap}}(S)$, then $X \cdot Y \in \mathsf{Pl_{Ap}}(S)$, where
\[
	X \cdot Y = \{xy \in S \mid x\in X, y\in Y\}.
\]
\end{enumerate}
Hence $\mathsf{Pl_{Ap}}(S)$ is a subsemigroup of the semigroup $P(S)$, consisting of all subsets of $S$ under the above multiplication (3).

Here is a summary of some further properties. For further details, see~\cite{qtheory}. For point (8), see~\cite{AHNR.1995}.

\bp\mbox{}
\label{proposition.properties}
\begin{enumerate}
\item (push) If $\varphi: S \to T$ is a surmorphism and $X \in \mathsf{Pl_{Ap}}(S)$, then $\varphi(X) \in \mathsf{Pl_{Ap}}(T)$.
\item (pull) If $\varphi: S \to T$ is a surmorphism and $Y\in \mathsf{Pl_{Ap}}(T)$, then there exists $X \in \mathsf{Pl_{Ap}}(S)$
such that $\varphi(X)=Y$.
\item (monadic) The map
\[
\begin{split}
	U: \quad \mathsf{Pl_{Ap}}(\mathsf{Pl_{Ap}}(S)) &\to \mathsf{Pl_{Ap}}(S)\\
	\{X_1,\ldots,X_k\} & \mapsto X_1 \cup \cdots \cup X_k
\end{split}
\]
is a surmorphism.
\item \label{injective ssg} (injective subsemigroups)
We have the following injective inclusions of semigroups
\[
	S \leqslant \mathsf{Pl_{Ap}}(S) \leqslant \mathsf{Pl_{Ap}}(\mathsf{Pl_{Ap}}(S)) \leqslant \cdots.
\]
\item (injections with respect to relational morphisms) If
\[
	r: S \to T \in \mathsf{V}
\]
is a relational morphism, then $r$ extends to $\mathsf{Pl_{Ap}}(S) \to T \in \mathsf{V}$ by sending $X \in \mathsf{Pl_{Ap}}(S)$ to $r(X)$.
\item \label{compact} (compactness) There exists a relational morphism $r: S \to T \in \mathsf{V}$ such that
$X \in \mathsf{Pl_{Ap}}(S)$ if and only if there exists a $T \in \mathsf{V}$ such that $X \subseteq r^{-1}(T)$. We say that
$r$ {\em computes} $\mathsf{Pl_{Ap}}(S)$.
\item (singletons) We have $S\in \mathsf{V}$ if and only if $\mathsf{Pl_{Ap}}(S) = \{\emptyset, \{s\} \mid s\in S\}$.
\item (decidability) There exists a pseudovariety $\mathsf{V}$ such that membership of finite semigroups in $\mathsf{V}$ is decidable, 
but the computation of $\mathsf{Pl_{Ap}}(S)$ for arbitrary S is not decidable.
\end{enumerate}
\ep

Note that, regarding Proposition~\ref{proposition.properties} \eqref{compact}, there exist examples (see Example~\ref{example.pointlike} below),
where $r$ computes $\mathsf{Pl_{Ap}}(S)$, but $r$ extended to $\mathsf{Pl_{Ap}}(S)$ does not compute $\mathsf{Pl_{Ap}}(\mathsf{Pl_{Ap}}(S))$. How far up the 
sequence in Proposition~\ref{proposition.properties}~\eqref{injective ssg} the induced map $r$ computes $\mathsf{Pl_{Ap}}^{(k)}(S)$ is
a measure for the {\em tightness} of $r$.

\be
\label{example.pointlike}
Let $S = Z_2 = \{1,-1\}$. Then the relational morphism $r: Z_2 \to 1$ computes
$\mathsf{Pl_{Ap}}(Z_2) = P(Z_2)$. The induced map $r$ is trivially $P(Z_2) \to \{0,1\}$, where the empty set is mapped to $0$.
This implies that $\{\emptyset, Z_2,1,-1\}$ is in $\mathsf{Pl_{Ap}}( \mathsf{Pl_{Ap}}(Z_2))$, if the 
induced map $r$ computed $\mathsf{Pl_{Ap}}( \mathsf{Pl_{Ap}}(Z_2))$. But 
$\{\emptyset, Z_2,1,-1\}$ maps onto an aperiodic semigroup by identifying $1$ and $-1$, so this set is not pointlike.
\ee

\begin{itemize}
\item{{\color{blue} The Karnofsky--Rhodes Expansions}   

An {\em expansion} is a functor $F$ on the category of finite semigroups such that there is a surjective morphism $F(S) \rightarrow S$ that forms part 
of a natural transformation from $F$ to the identity functor. Expansions have played a crucial role in finite semigroup theory, the first being the 
Rhodes expansion, that was used for a proof of the Fundamental Lemma of Complexity. See \cite{TilsonXII}. Expansions intuitively remove ``singularities'' in semigroups 
replacing them with an expanded semigroup with better behavior. In the proof for computability of complexity, expansions are used to find a cover of 
an aperiodic semigroup whose Cayley graph look like the Cayley graph of free aperiodic Burnside semigroups. We will discuss this below.

The first expansion we look at is the Karnofsky--Rhodes expansion. We briefly recall its definition. See~\cite[Definition 4.15]{gst}. 
Let $S$ be a semigroup with a generating set $X$. The right Cayley graph $\operatorname{RCay}(S,X)$ of $S$ with respect to X is the rooted 
graph with vertex set $S^\mathcal{I}$, where $S^\mathcal{I}$ adjoins an identity $\mathcal{I}$ to $S$, even if $S$ has an identity. It has root 
$r=\mathcal{I}$ and edges $s\stackrel{a}{\rightarrow}sa$ for all $(s, a, s^{\prime}) \in S^\mathcal{I} \times X \times S^\mathcal{I}$, where 
$s^{\prime} = sa$ in $S^\mathcal{I}$. Notice that the strongly connected components of the right Cayley graph have vertices the 
$\mathcal{R}$-classes of $S^{\mathcal{I}}$. An edge $s\stackrel{a}{\rightarrow}sa$ is called a {\em transition edge} if 
$sa <_{\mathcal{R}} s$ in $S^{\mathcal{I}}$. That is, the edge goes from one strongly connected component to one lower in the usual order on connected 
components of a directed graph. The {\em right Karnofsky--Rhodes expansion} identifies two non-empty paths starting at $\mathcal{I}$ if their 
edges multiply to the same element of $S$ and traverse the same set of transition edges. This defines a congruence on the free semigroup $X^{+}$. 
The quotient semigroup $\operatorname{RKR}(S,X)$ is called the {\em right Karnofsky--Rhodes expansion relative to $X$}. There is an obvious 
dual construction of the left Karnofsky--Rhodes expansion based on the left Cayley Graph that we denote by $\operatorname{LKR}(S,X)$.}

We will need a version of the Karnofsky--Rhodes expansion for more general automata than the Cayley automaton of a semigroup $S$. Assume that we have a complete deterministic accessible automaton $\mathcal{A}=(Q,X,\mathcal{I})$ with start state $\mathcal{I}$. Furthermore assume that no transition in the automaton ends at $\mathcal{I}$. Let $S$ be the transition semigroup of $\mathcal{A}$ so that we have the corresponding ts $(Q,X)$ that has a distinguished state $\mathcal{I}$ such that for all $q \in Q$ there is an $s \in S^{\mathcal{I}}$ such that $\mathcal{\mathcal{I}}s=q$. By assumption, $\mathcal{I}$ is the unique element $s$ of $S^{\mathcal{I}}$ such that $\mathcal{I}s=\mathcal{I}$.

There is an obvious notion of strongly connected component and of transition edges in $\mathcal{A}$. The {\em  Karnofsky--Rhodes expansion} of $\mathcal{I}$ identifies two words $w,v$ in $X^{+}$ if the path they read starting at $\mathcal{I}$ satisfies $\mathcal{I}v=\mathcal{I}w$ 
 and they traverse the same set of transition edges. This defines a right congruence $\kappa$ on the free semigroup $X^{+}$. 
The quotient semigroup $\operatorname{KR}(\mathcal{A})$ is called the {\em  Karnofsky--Rhodes expansion of $\mathcal{A}$}. We also have the quotient ts $(X^{+}/\kappa,T)$ where $T$ is the image of $X^+$ under its action on $X^{+}/\kappa$. The Karnofsky-Rhodes expansion is indeed an expansion in the category of ts considered above.
\end{itemize}

We give an example of the Karnofsky-Rhodes expansion.

\be

\label{example.KR}
Consider the right Cayley graph of the Klein $4$-group $Z_2 \times Z_2$ with zero with generators 
$\{a,b,\square\}$, where  $a=(1,-1)$, $b=(-1,1)$, and $\square$ is the zero. The right Cayley graph
$\mathsf{RCay}(Z_2 \times Z_2 \cup \{\square\},\{a,b,\square\})$ is
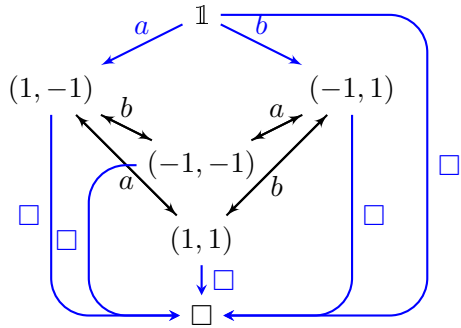
\begin{center}
\begin{tikzpicture}[auto]
\node (A) at (0, 0) {$\bf{1}$};
\node (B) at (-2,-1) {$(1,-1)$};
\node(C) at (2,-1) {$(-1,1)$};
\node(D) at (0,-2) {$(-1,-1)$};
\node(E) at (0,-3) {$(1,1)$};
\node(F) at (0,-4) {$\square$};
\draw[edge,blue,thick] (A) -- (B) node[midway, above] {$a$};
\draw[edge,blue,thick] (A) -- (C) node[midway, above] {$b$};
\draw[edge,thick] (B) -- (D) node[midway, above] {$b$};
\draw[edge,thick] (D) -- (B);
\draw[edge,thick] (C) -- (D) node[midway, above] {$a$};
\draw[edge,thick] (D) -- (C);
\draw[edge,thick] (B) -- (E) node[midway, below] {$a$};
\draw[edge,thick] (E) -- (B);
\draw[edge,thick] (C) -- (E) node[midway, below] {$b$};
\draw[edge,thick] (E) -- (C);
\draw[arrow, rounded corners=5mm,blue] (A) -- ($(A) + (3,0)$) -- ($(C)+(1,0)$)  node[pos=2, right]{$\square$}
-- ($(F) + (3,0)$)  -- (F); 
\draw[arrow, rounded corners=5mm,blue] (C) -- ($(F) + (2,0)$) node[pos=0.5, right]{$\square$} -- (F); 
\draw[arrow, rounded corners=5mm,blue] (B) -- ($(F) + (-2,0)$) node[pos=0.5, left]{$\square$} -- (F); 
\draw[arrow, rounded corners=5mm,blue] (D) -- ($(D) + (-1.5,0)$) -- ($(F) +(-1.5,0)$) node[pos=0.5, left]{$\square$} -- (F) ; 
\draw[arrow,blue] (E) -- (F) node[midway,right]{$\square$}; 
\end{tikzpicture}
\end{center}
where all three arrows $a,b,\square$ fix the vertex $\square$ at the bottom.
Transition edges are indicated in blue. Double edges mean that right multiplication by the label for either vertex yields 
the other vertex. The right Karnofsky--Rhodes expansion of this right Cayley graph is given by
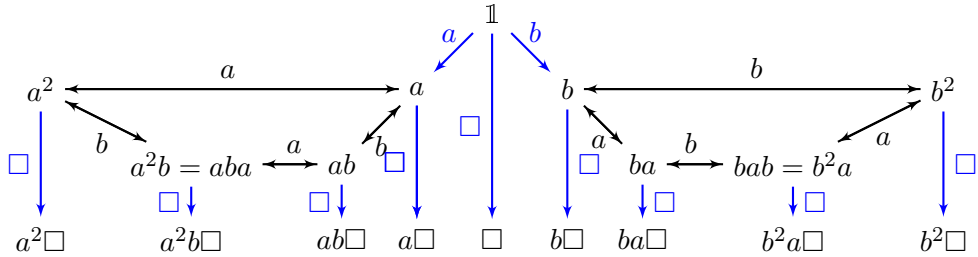
\begin{center}
\begin{tikzpicture}[auto]
\node (A) at (0, 0) {$\bf{1}$};
\node (B) at (-1,-1) {$a$};
\node(C) at (1,-1) {$b$};
\node(D) at (-2,-2) {$ab$};
\node(E) at (2,-2) {$ba$};
\node(F) at (-6,-1) {$a^2$};
\node(G) at (6,-1) {$b^2$};
\node(H) at (-4,-2) {$a^2b=aba$};
\node(I) at (4,-2) {$bab=b^2a$};
\node(ZA) at (0,-3) {$\square$};
\node(ZB) at (-1,-3) {$a\square$};
\node(ZD) at (-2,-3) {$ab\square$};
\node(ZH) at (-4,-3) {$a^2b\square$};
\node(ZF) at (-6,-3) {$a^2\square$};
\node(ZC) at (1,-3) {$b\square$};
\node(ZE) at (2,-3) {$ba\square$};
\node(ZI) at (4,-3) {$b^2a\square$};
\node(ZG) at (6,-3) {$b^2\square$};
\draw[edge,blue,thick] (A) -- (B) node[midway, above] {$a$\;};
\draw[edge,thick,blue] (A) -- (C) node[midway, above] {\;$b$};
\draw[edge,thick] (B) -- (F) node[midway,above] {$a$\;};
\draw[edge,thick] (F) -- (B);
\draw[edge,thick] (B) -- (D) node[midway,below] {$b$};
\draw[edge,thick] (D) -- (B);
\draw[edge,thick] (F) -- (H) node[midway,below] {$b$\;};
\draw[edge,thick] (H) -- (F);
\draw[edge,thick] (D) -- (H) node[midway,above] {\;$a$};
\draw[edge,thick] (H) -- (D);
\draw[edge,thick] (C) -- (E) node[midway,below] {$a$\;};
\draw[edge,thick] (E) -- (C);
\draw[edge,thick] (E) -- (I) node[midway,above] {$b$\;};
\draw[edge,thick] (I) -- (E);
\draw[edge,thick] (C) -- (G) node[midway,above] {\;$b$};
\draw[edge,thick] (G) -- (C);
\draw[edge,thick] (G) -- (I) node[midway,below] {\;$a$};
\draw[edge,thick] (I) -- (G);
\draw[edge,blue,thick] (A) -- (ZA) node[midway, left] {$\square$};
\draw[edge,blue,thick] (B) -- (ZB) node[midway, left] {$\square$};
\draw[edge,blue,thick] (D) -- (ZD) node[midway, left] {$\square$};
\draw[edge,blue,thick] (H) -- (ZH) node[midway, left] {$\square$};
\draw[edge,blue,thick] (B) -- (ZB) node[midway, left] {$\square$};
\draw[edge,blue,thick] (F) -- (ZF) node[midway, left] {$\square$};
\draw[edge,blue,thick] (C) -- (ZC) node[midway, right] {$\square$};
\draw[edge,blue,thick] (E) -- (ZE) node[midway, right] {$\square$};
\draw[edge,blue,thick] (I) -- (ZI) node[midway, right] {$\square$};
\draw[edge,blue,thick] (G) -- (ZG) node[midway, right] {$\square$};
\end{tikzpicture}
\end{center}
where arrows $a,b,\square$ fix all the vertices at the bottom.
\ee

\begin{itemize}
\item{{\color{blue} The McCammond Expansion and the Geometric Semigroup Theory (gst) Expansion}

We will define  the McCammond expansion. For more details on the McCammond expansion,
 see~\cite{gst}. Recall that a simple path in a graph or an automaton is a path that does not visit any vertex twice. Empty paths are considered simple. 
 Let $\mathcal{A}$ be an automaton. We assume that $\mathcal{A}$ is deterministic and has a unique start state $I$. The McCammond expansion 
 $\operatorname{Mc(\mathcal{A})}$ of $\mathcal{A}$ 
 is the automaton with vertex set $V$, which is the set of simple paths in $\mathcal{A}$. The edges are given by the following rules. If 
 $p$ is a simple path and $a \in X$ is such that $p(p,a,pa)$ is a simple path, then we have an edge from $p$ to $p(p,a,pa)$. Otherwise, there is a 
 unique vertex $v$ on $p$ such that $v=pa$. In this case, letting $q$ be the prefix of $p$ from its start vertex to $v$ we have an edge $(p,a,q)$. 
 Its start state is the empty path at $I$.

It is known that $\operatorname{Mc}(\mathcal{A})$ has the unique simple path property, defined as follows:  A rooted graph $(\Gamma, r)$ 
with root $r$ has the unique simple path property if for each vertex $v$ in $(\Gamma,r)$ there is a unique simple path from the root $r$ to $v$. 
The transition semigroup of $\operatorname{Mc}(\mathcal{A})$ is known to be aperiodic if $S$ is aperiodic~\cite{gst}.

The Geometric Semigroup Theory (gst) expansion of $(S,X)$ is the composition of the McCammond 
expansion and the right Karnofsky--Rhodes expansion. That is,
\[
	\operatorname{gst}(S,X) = \operatorname{Mc}\circ \operatorname{RKR}(S,X).
\]
The $\operatorname{gst}$ expansion is used extensively in this paper and in \cite{complexityn}.
}
\end{itemize}

The McCammond expansion is not functorial on the category of (transformation) semigroups~\cite{gst}. Despite this we show how to extend a 
relational morphism $r\colon S \rightarrow T$ of semigroups to a relational morphism 
$\tilde{r}\colon \operatorname{Mc}(S) \rightarrow \operatorname{Mc}(T)$. An obvious analogous construction works in the case of transformation semigroups.

\bd \label{Mcexp_rm}

Let $r\colon S \rightarrow T$ be a relational morphism, and $\alpha\colon\operatorname{Mc}(S) \rightarrow S$ and 
$\beta\colon\operatorname{Mc}(T) \rightarrow T$ be the McCammond covering maps of $S$ and $T$. The {\em McCammond expansion} 
of $r$ is the relational morphism $\tilde{r}=\alpha r \beta^{-1} \colon \operatorname{Mc}(S) \rightarrow \operatorname{Mc}(T)$.

\ed
\noindent
See \cite{RS1} where both the Karnofsky--Rhodes expansion and the McCammond expansion are used in the theory of Markov Chains.

\begin{itemize}
\item{{\color{blue} Geometric properties of Cayley graphs of Free Aperiodic Semigroups}  

The Free Semigroup $F_{n}(X)$ generated by a set $X$ and index $n$ is the free semigroup in the variety defined by the identity $x^{n}=x^{n+1}$ 
and generated by $X$.

Starting with work of McCammond in 1991 \cite{Mc91}, these and other free Burnside semigroups have been extensively studied. 
For $n>2$ their word problems are decidable.

Most important for the proof of decidability of the complexity is that for $n>2$ the semigroup $F_{n}(X)$ is finite $\mathcal{J}$-above.} 
\end{itemize}

\begin{itemize}
\item{{\color{blue} Methods from Profinite Semigroup Theory} 

Profinite semigroup theory has played a central role in finite semigroups since the 1980s. See~\cite{Almeida:book} for background and applications. 
In particular, there is a proof of the type II conjecture using new properties of free profinite groups. We need further profinite results concerning 
the notion of inevitable graphs that we build into the automata and semigroups we consider in the proof of computability of complexity in this paper and in~\cite{complexityn}. See~\cite{stablepairs} for the definition and basic properties of inevitable graphs.  
}
\end{itemize}

\noindent
An important application of profinite methods in complexity theory is in the theory of stable pairs \cite{Henckellstable}. See \cite{stablepairs} for a 
profinite proof of the main theorem of \cite{Henckellstable}.

\subsection{The Monoid of Closure Operations}\label{Closops}

We gather definitions and results from \cite{Trans}. For more details the reader is strongly urged to consult this paper. Since all the computations we do on examples in this document are done in the Evaluation Transformation Semigroup defined in~\cite[Section 5]{Trans} our goal is to summarize the background material in that paper needed to define this object. We begin with the definition of the monoid of closure operations on the direct product $L^{2} = L \times L$ of a lattice $L$ with itself. 

Let $L$ be a lattice and let $L^{2} = L \times L$. Let $f$ be a closure operator on $L^2$. By definition this means that $f$ is an order preserving, extensive (that is, for all $(l_{1},l_{2}) \in L^{2}, (l_{1},l_{2}) \leqslant (l_{1},l_{2})f$), idempotent function on $L^2$. A {\em stable pair} for $f$ is a closed element of $f$. Thus a stable pair $(l,l') \in L^2$ is an element such that $(l,l')f = (l,l')$.  The stable pairs of $f$ are a meet closed subset of $L \times L$. Conversely, each meet closed subset of $L \times L$ is the set of stable pairs for a unique closure operator on $L \times L$. We will identify $f$ as a binary relation on $L$ whose pairs are precisely the stable pairs. Let $B(L)$ be the monoid of binary relations on the set $L$.  As is well known, $B(L)$ is isomorphic to the monoid $ M_{n}(\mathcal{B})$ of $n \times n$ matrices over the 2-element Boolean algebra $\mathcal{B}$, where $n=|L|$. The Boolean matrix associated to $f$ is of dimension $|L| \times |L|$. It has a 1 in position $(l_{1},l_{2})$ if $(l_{1},l_{2})$ is a stable pair and a 0 otherwise. 

With this identification, the collection $\mathcal{C}(L^{2})$ of all closure operators on $L^2$ is a submonoid of the monoid $B(L)$ of binary relations on $L$~\cite[Proposition 2.5]{Trans}. We thus also consider $\mathcal{C}(L^{2})$ to be a monoid of $|L| \times |L|$ Boolean matrices. Let $L$ be either $Rh_{B}(G)$ or $\operatorname{SP}(G\times B)$. We now recall the definition of some important unary operations on $\mathcal{C}(L^{2})$. 
\bd

Let Let $f \in \mathcal{C}(L^{2})$.

\begin{enumerate}

\item{The domain of $f$ denoted by $\operatorname{Dom}(f)$ is the set $\{x \mid \exists y, (x,y) \in f\}.$}

\item{Define the relation $\overleftarrow{f}$ by $\overleftarrow{f}= \{(x,x) \mid x \in \operatorname{Dom}(f)\}$. $\overleftarrow{f}$ is called {\em back flow along $f$}. See \cite[Remark 2.26]{Trans} for the reason for this terminology.}

\item{The relation $f^{*}$ is defined by $f^{*} = f \cap \{(x,x) \mid x \in X\}$. $f^{*}$ is called {\em the Kleene closure of $R$}. See~\cite[Sections 2 and 4]{Trans} for the reason for this terminology.}

\item{Define the {\em loop of $f$} to be the relation $f^{\omega+*}=f^{\omega}f^{*}$, where $f^{\omega}$ is the unique idempotent in the subsemigroup generated by $f$.}

\end{enumerate}

\ed

At times it is convenient to identify $\overleftarrow{f}$ as $1|_{\operatorname{Dom}(f)}:L \rightarrow L$, the identity function restricted to 
$\operatorname{Dom}(f)$. Similarly, we identify $f^{*}$ as $1|_{\operatorname{Fix{f}}}:L \rightarrow L$, the restriction of the identity to the set of fixed-points of $f$, where $x \in L$ is a fixed point if $(x,x) \in f$. Our use of these will be clear from the context.

\subsection{The 0-Flow Monoid}\label{FMon}

Let $S$ be a $\GM$ semigroup generated by $X$ and let $x \in X$. Let $I(S) = \mathcal{M}^{0}(A,G,B,C)$. We define a binary relation $f_{x}$
on $\operatorname{SP}(G \times B)$ by $((Y,\Pi), (Z,\Theta)) \in f_{x}$ if and only if $Yx \subseteq Z$ and the partial function induced by
right multiplication by $x$, $\cdot x: Y \rightarrow Z$ induces a well-defined partial injective map $\cdot x: Y/\Pi \rightarrow Z/\Theta$. This means that if  $(g,b), (g',b') \in Y$ and $(g,b)x,(g',b')x$ are both defined (and hence in $Z$ ), then $(g,b)\Pi (g',b')$ if and only if $(g,b)x\Theta (g',b')x$. Then $f_{x} \in \mathcal{C}(L^{2})$. See \cite[Proposition 2.22]{Trans}. $f_{x}$ is called the free-flow by $x$.

Let $V$ be a pseudovariety of monoids. We need to recall the definition of a $V$-stable pair in a finite monoid.

\bd
Let $V$ be a pseudovariety and let $S$ be a monoid. Suppose that $A \subseteq S$ and
$S'$ is a submonoid of $S$. Then $(A,S')$ is called a $V$-stable pair if, for all relational
morphisms $\phi\colon S \rightarrow T$ with $T \in V$, there is an element $t \in T$ so that $A\subseteq t\phi^{-1}$ and
$S' \subseteq \operatorname{Stab}(t)$, where $\operatorname{Stab}(t)=\{u \in T \mid tu=t\}$ is the right stabilizer of $t$ in $T$.
\ed

It is proved in \cite{stablepairs, Henckellstable} that if $Ap$ is the variety of aperiodic monoids, then for all finite monoids $S$ and all pairs $(A,S')$, where $A \subseteq S$ and $S'$ is a submonoid of $S$ it is decidable if $(A,S')$ is an $Ap$-stable pair.

We now define the 0-flow monoid $M_{0}(L)$ as follows. 

\bd\label{Flowops}

Let $L = \operatorname{SP}(G \times B)$. The 0-flow monoid $M_{0}(L)$, is the
smallest subset of $\mathcal{C}(L^{2})$ satisfying the following axioms:

\begin{enumerate}

\item {(Identity) The multiplicative identity $I$ of $\mathcal{C}(L^{2})$ is in $M_{0}(L)$.}

\item{(Points) For all $x \in X$, $f_{x}$ the free-flow along $x$ belongs to $M_{0}(L)$.}

\item{(Products) If $f_{1}, f_{2} \in M_{0}(L)$, then $f_{1}f_{2} \in M_{0}(L)$.}

\item{(Vacuum) If $f \in M_{0}(L)$, then $\overleftarrow{f} \in M_{0}(L)$.}\label{Vac}

\item{(Loops) If $f \in M_{0}(L)$, then $f^{\omega+*}\in M_{0}(L)$.}

\end{enumerate}
\ed

We remark that for each $n \geqslant 0$, there is an $n$-flow monoid $M_{n}(L)$ defined in \cite{Trans}. The definition above of $M_{0}(L)$ is exactly what is defined in \cite{Trans}. For $n>0$, Axioms (1)-(4) are the same for $M_{n}(L)$ as for $M_{0}(L)$. Axiom (5) restricts the use of the loop operator to $n$-loopable elements \cite[Section 4]{Trans}. In \cite{complexityn}, $n$-loopable elements are replaced by a more restrictive definition and the modified $M_{n}(L)$ plays a crucial role in the main results of \cite{complexityn}.

%%%%%%%%%%%%%%%%%%%%%%%%%%%%%%%%%%%%%%%%%%%%%%%%%%%%%%%
\subsection{The Evaluation Transformation Semigroup}\label{Ets}

We now define the Evaluation Transformation Semigroup $\mathcal{E}(L) = (\States, \Eval(L))$ as in \cite{Trans}. Again, because of our interest in this paper on aperiodic flows, we do not define the $n$-Evaluation Transformation Semigroups $E_{n}$ for $n>0$. We begin with the definition of Well-Formed Formulae (WFFs).

\bd

Let $X$ be an alphabet. We define a well-formed formula inductively as follows.

\begin{enumerate}

\item{The empty string $\epsilon$ is a well-formed formula.}

\item{Each letter $x \in X$ is a well-formed formula.}

\item{If $\tau, \sigma$ are well-formed formulae, then
so is $\tau\sigma$.}

\item{If $\tau$ is a well-formed formula that is not a proper power (i.e., not of the form
$\sigma^{n}$ where $n > 1$), then $\tau^{\omega+*}$ is also a well-formed formula.} 

\end{enumerate}

\ed

The set of well-formed
formulae is denoted by $\Omega(X)$. Well-formed formulae will be denoted by Greek
letters. As a convention, if $\tau=\sigma
^{n}$, where $\sigma$ is not a proper power, then we set
$\tau^{\omega+*}=\sigma^{\omega+*}$. In other words, we extract roots before applying the unary operation $\omega+*$.

Let $V = \prod_{f \in M_{0}(L)}\overleftarrow{f}$. $V$ is called the {\em Vacuum}. See \cite{Trans} for an explanation of this terminology. In \cite{Trans}, $V$ is denoted by $\mathcal{F}_{0}$. It is proved in \cite{Trans} that $V$ is an idempotent in $M_{0}(L)$. We will now work exclusively in the subsemigroup $VM_{0}(L)V$ of $M_{0}(L)$. We want to interpret WFFs in $VM_{0}(L)V$.

\bd

Define recursively a partial function $\mathcal{I}: \Omega(X)\rightarrow VM_{0}(L)V$ as follows. 

\begin{enumerate}

\item{$\epsilon\mathcal{I} = V$.}

\item{$x\mathcal{I} = VxV$ for $x \in X$.}

\item{If $\mathcal{I}$ is already defined on $\tau, \sigma \in \Omega(X)$, set $(\tau\sigma)\mathcal{I} = \tau\mathcal{I}\sigma\mathcal{I}$.}

\item{If $\tau \in \Omega(X)$ is not a proper power and $\tau\mathcal{I}$ is defined, set $\tau^{\omega+*}\mathcal{I} = 
(\tau\mathcal{I})^{\omega+*}$.}

\end{enumerate}

\ed

We normally omit $\mathcal{I}$ and assume that a WFF $\tau$  is being evaluated in $VM_{0}(L)V$ according to the definition of $\mathcal{I}$. 
 We first define a new operator on elements of $M_{0}(L)$ called {\em forward flow}. Recall that the bottom of the lattice $\operatorname{SP}(G \times B)$ is the pair $\Box = (\emptyset, \emptyset)$.

\bd

Let $f \in M_{0}(L)$ and let $l \in L$. Let $(l,\Box)f =(l_{1},l_{2})$ We define the forward flow of $f$ denoted by $\overrightarrow{f}$ by 
$l\overrightarrow{f}=l_{2}$. That is, we apply $f$ to $(l,\Box)$ and project to the right-hand coordinate.

\ed

The following is proved in~\cite[Sections 2 and 4]{Trans}.

\begin{enumerate}

\item{$\overrightarrow{f}:L \rightarrow L$ is an order preserving function on $L$.}

\item{If $lV=l$, then $(l,\Box)f=(l,l')$ and $l' \in LV$. That is, the left-coordinate of $(l,\Box)f$ is still $l$ and the right-hand coordinate is in $LV$. Therefore $\overrightarrow{f}$ is a well-defined function from $LV$ to $LV$.}

\item{The assignment of $f$ to $\overrightarrow{f}$ defines an action of $VM_{0}(L)V$ on $LV$. It follows that we have a transformation semigroup $(LV,M'_{0}(L))$, where $M'_{0}(L)$ is the image of $VM_{0}(L)V$ on $LV$ under this action.}

\end{enumerate}

We now restrict the action in $(LV,M'_{0}(L))$ to the set of $\States$ defined as follows. Let $(g,b) \in G \times B$. Then 
the element $(\{(g,b)\},\{(g,b)\}) \in L$ is called a {\em point}. In~\cite[Section 5]{Trans}, it is proved that every point $p$ satisfies $pV=p$. 

\bd

The set of $\States$ is the smallest subset of $LV$ such that:

\begin{enumerate}

\item{Every point $p \in \States$.}

\item{If $l \in \States$, then $l\overrightarrow{f} \in \States$.}

\end{enumerate}

\ed

In other words $\States$ is the smallest subset of $LV$ containing the points and closed under the action of $M_{0}(L)$ on $LV$. We can finally define the Evaluation Transformation Semigroup where all the computations in the examples in this paper take place.

\bd

The Evaluation Transformation Semigroup $\mathcal{E}(L)$ is defined by $\mathcal{E}(L) = (\States, \Eval(L))$ where $\Eval(L)$ is the image of $M'_{0}(L)$ by restricting its action to $\States$.

\ed

We remark that there is an Evaluation Transformation Semigroup $\mathcal{E}_{n}(L)$ for all $n \geqslant 0$. The definition here is the case $n=0$
which is all we need in this paper as we are only concerned with aperiodic flows.   
%---------------------------------------------------------------------------------------------------
\section{Outline of the Proof for Computability of Complexity 1}
\label{section.outline}
%---------------------------------------------------------------------------------------------------

Let $(G \times B,S)$ be $\GM$ and let $X$ be a generating set for $S$. 
Our main theorem Theorem \ref{MainTh} can be rephrased as follows. For complexity 1, it says that the lower bound defined in~\cite[Section 5]{Trans} is perfect.

\bt
\label{cequal1}
Let $(G \times B,S)$ be $\GM$. Then $Sc =1$ if and only if $\RLM(S)c \leqslant 1$ and every element of $\States$ is a cross-section.
\et

\noindent
Here are the main steps in the proof.
\begin{enumerate}
\item{Based on $(\States,\Eval(L))$ we effectively define an aperiodic automaton $\mathcal{A}$ with state set $Q$. This definition involves 
all the tools we described in the previous section.}
\item{A partial flow is a partial function $F \colon Q \rightarrow \operatorname{SP}(G \times B)$ that satisfies the properties in the definition of a 
flow for all elements in $\operatorname{Dom}(F)$. We define iteratively a sequence $F_{i}$ for $i=0,1,\ldots$ of  partial flows with Domains 
$D_{i} \subseteq Q$ with $D_{i} \subseteq D_{i+1}$ that approximate a flow from $\mathcal{A}$ to the state set of $\Eval(L)$.}
\item{We prove that our iterative procedure either converges to a flow in a computable finite number of steps from our aperiodic automaton to 
$\operatorname{SP}(G \times B)$ if every element in $\States$ is a cross-section or there is no automaton with aperiodic transformation semigroup 
that defines a flow.}
\end{enumerate}

%%%%%%%%%%%%%%%%%%%%%%%%%%%%%%%%%%%%%%%%%%%%%%%%%%%%%
\section{Profinite smoothing}  
\label{section.6}
%%%%%%%%%%%%%%%%%%%%%%%%%%%%%%%%%%%%%%%%%%%%%%%%%%%%%

The goal of this section is to prove Theorem \ref{MainTh} and thus computing that a semigroup has complexity 1 is decidable.

%%%%%%%%%%%%%%%%%%%%%%%%%%%%%%%%%%%%%%%%%%%%%%%%%%%%%%%%%%%%
\subsection{Background in free aperiodic Burnside semigroups and free aperiodic $\omega$-semigroups}\label{freeaperiodic}

We need some background in free aperiodic Burnside semigroups and free aperiodic $\omega$-semigroups. We review the material we need in this paper. In~\cite{Mc91} McCammond treats free aperiodic Burnside semigroups. In~\cite{Mc-nf} he treats free aperiodic $\omega$-semigroups and gives the connection to free aperiodic Burnside semigroups. We begin by reviewing free aperiodic Burnside semigroups.  

Let $n$ and $m$ be positive integers. Let $X$ be a set of cardinality $n$. We define $B_{n}(m)$, the free Burnside $n$-generated aperiodic semigroup 
of degree $m$ to be the free semigroup with $n$-generators in the variety of semigroups defined by the identity $x^{m}=x^{m+1}$. Up to isomorphism, 
$B_{n}(m)$ is the quotient of the free semigroup $X^+$ by the smallest congruence containing 
$\{(w^{m},w^{m+1}) \mid w \in X^{+}\}$. We will always assume that we have a fixed $X$ generated $\GM$ semigroup $(S,X)$, so that $n$ is fixed 
throughout the discussion. The variable $m$ will always be an integer $m \geqslant 6$ and will vary in the discussion depending on the context.

Let $\mathsf{Cay}_{m}(X)$ be the right Cayley graph of $B_{n}(m)$ relative to the generating set $X$. We are working with semigroup Cayley graphs. So $\mathsf{Cay}_{m}(X)$ has vertices the elements of $B_{n}(m)$ 
and an initial state $I$ not belonging to $B_{n}(m)$. We  identify $x$ with its value in $B_{n}(m)$ under the canonical morphism from $X^{+}$. 
The ($X$ labeled) edges of $\Gamma$ are of the form $s \stackrel{x}{\rightarrow} sx$, where for $s = I$, we define $Ix$ to be $x$. Notice that 
this implies that there are no edges that end at $I$.

An effective construction of $\mathsf{Cay}_{m}(X)$ solves the word problem for $B_{n}(m)$. We do this by considering for each $s \in B_{n}(m)$, the subautomaton 
$\mathsf{str}_m(s)$ of $\mathsf{Cay}_{m}(X)$ (called the straightline automaton for $s$ in  $\mathsf{Cay}_{m}(X)$) that has $I$ as initial state and $s$ in $B_{n}(m)$ as terminal state and all states that lie on a path from $I$ to $s$. Note that the states of $\mathsf{str}_m(s)$ are those 
in the quotient of $\mathsf{Cay}_{m}(X)$ by the right ideal of elements $t \in B_{n}(m)$ such that $t <_{\mathcal{R}} s$. Thus the states of $\mathsf{str}_m(s)$ 
consist of $I$, all elements $t$ such that $s \leqslant_{\mathcal{R}} t$ and a sink state $\square$ representing all the elements in the right ideal above.

One of the very important results of~\cite{Mc91} is that for every $s \in B_{n}(m)$, the automaton $\mathsf{str}_m(s)$ is a finite automaton. Moreover,
$\mathsf{str}_m(s)$ accepts the word problem for $s$, that is the language $WP(s) = \{w \in X^{+} \mid w\theta = s\}$, where 
$\theta \colon X^{+} \rightarrow B_{n}(m)$ is the natural morphism. It follows that $WP(s)$ is a computable regular language for every $s \in B_{n}(m)$. 
McCammond also constructs a regular expression for $WP(s)$ that does not use the union operation. 
For each word $w$, there is a unique minimal length word $\mathsf{red}(w)$ such that $w\theta = \mathsf{red}(w)\theta \in B_n(m)$. Furthermore if 
$\mathsf{red}(w) = x_1x_2 \cdots x_k$ for $x_i \in X$, then it is shown in~\cite{Mc91,Mc-nf} that every directed loop in  $\mathsf{str}_m(w\theta)$
not containing $\square$ passes through at least one of the vertices $(x_1\cdots x_j)$ for  $1 \leqslant  j \leqslant  k$. By definition, no non-empty loop passes through $I$. In the language of \cite{gst}, $\mathsf{str}_{m}(s)$ is a linear automaton for every $s \in B_{n}(m)$. This means that the natural order on its strongly connected components is a linear order and that the transition edges (edges between strongly connected components) form a chain.

McCammond~\cite{Mc91,Mc-nf} constructs this automaton in two ways. The first 
is a detailed study of the Todd--Coxeter process that gives a sequence of approximating finite automata for $\mathsf{str}_m(s)$ whose limit 
(in an appropriate sense) is all of $\mathsf{str}_m(s)$. The difficult part of the theorem is to prove that this process ends in a finite number of steps. 
After doing this a Church--Rosser rewriting system is constructed for the word problem. The system computes $\mathsf{red}(w)$ for each $w \in X^+$. It is shown how to compute both the rewriting system and how to use it to effectively construct $\mathsf{str}_m(s)$ for all 
$s \in B_{n}(m)$. It is proved that there is a unique simple path in $\mathsf{str}_{m}(s)$ from $I$ to $s$ and that the normal form computed by the Church--Rosser system is precisely the word read by this path, which is the unique shortest word accepted by $\mathsf{str}_m(s)$.

We turn to the free aperiodic $\omega$-semigroup $\mathsf{\Omega}_{Ap}(X)$. By definition this is the smallest subsemigroup of the free profinite aperiodic semigroup generated $\widehat{F}_{Ap}(X)$ generated by $X$ that contains $X$ and is closed under the operation that sends $s \in \widehat{F}_{Ap}(X)$ to 
$s^{\omega}$, the unique idempotent in the closed subsemigroup generated by $s$.

In $\cite{Mc-nf}$ McCammond solves the word problem for $\mathsf{\Omega}_{Ap}(X)$ as a quotient of the free unary semigroup $U(X)$ by producing normal forms. For ease of reading we adopt McCammond's notation for elements of $U(X)$ that writes $(s)$ for $s^{\omega}$. This allows us to consider elements of $U(X)$ to be words over the alphabet $X \cup \{(,)\}$ with the usual compatibility condition for open and closed parentheses. McCammond's uses the following rules to transform an arbitrary element of $U(X)$ to its normal form by applying the following rules a finite number of times to elements $w,v \in U(X)$:
\begin{enumerate}
\item{$((w))=(w)$},
\item{$(w^{k})=(w)$ for all $k>0$},
\item{$(w)(w)=(w)$},
\item{$(w)w=(w)=w(w)$},
\item{$(wv)w=w(vw)$}.
\end{enumerate}

The specifics of the normal form are described in \cite{Mc-nf}. We do not need these details in the paper. We now note that $B_{n}(m)$ satisfies 
each of the rules 1-5 above when we interpret $(s)$ as $s^{m}$. Therefore there is a surjective generator preserving morphism 
$\theta_{m}\colon \mathsf{\Omega}_{Ap}(X) \rightarrow B_{n}(m)$. The following is~\cite[Theorem 13.1]{Mc-nf} and we use it crucially in this paper. 
It shows that $\mathsf{\Omega}_{Ap}(X)$ is a subdirect product of $B_{n}(m)$ as $m$ runs over the positive integers.

\bt\label{mainMcCam}

Let $s,t \in \mathsf{\Omega}_{Ap}(X)$ with $s \neq t$. Then there exists a computable integer $m$ such that $s\theta_{m} \neq t\theta_{m}$. 

\et

%%%%%%%%%%%%%%%%%%%%%%%%%%%%%%%%%%%%%%%%%%%%%%%%%%%%%%%%%%%%%%%%%%%

%%%%%%%%%%%%%%%%%%%%%%%%%%%%%%%%%%%%%%%%%%%%%%%%%%%%%%%

%%%%%%%%%%%%%%%%%%%%%%%%%%%%%%%%%%%%%%%%%%%%%%%%%%%%%%%
\subsection{Completion of automata and relational morphisms of automata}

From hereon in, complete automata always have a sink state $\square$, that is, a state such that $\square x =\square$ for all letters $x$. The sink state will be part of the signature of the automaton. Let $X$ be a finite set. The completion of a finite deterministic automaton $A=(I,Q,X,\cdot)$ with initial state $I$, state set $Q$ alphabet $X$ and action a partial function $\cdot\colon Q \times X \rightarrow Q$ is the complete automaton $A^{c}=(I,Q \cup\{\square\},X,\square)$ where $\square \notin Q$ is the sink state. If  $q \in Q, x \in X$ we have $qx=q\cdot x$ if this is defined in $A$, $qx=\square$ if $q\cdot x$ is not defined in $A$ and $\square x=\square$ for all $x \in X$. We emphasize that we add the sink state even if the action function of $A$ is a total function. 

We define the category $\mathcal{A}(X)$ to have as objects all isomorphism classes of finite complete trim deterministic 
automata $(I, Q, X, \square)$ with alphabet $X$. 
Trim means that 
there is a path in the automaton from $I$ to every state $q \in Q$. 
A morphism between $(I_1, Q_{1}, X, \square_{1})$ to $(I_2, Q_{2},X, \square_{2})$ is a function 
$f \colon Q_{1} \rightarrow Q_{2}$ such that $I_{1}f=I_{2}$ and for all $q \in Q, x \in X$ we have $(qx)f =(qf)x$ and $\square_{1}f=\square_{2}$. We write $\mathcal{A}$ for $\mathcal{A}(X)$ if the alphabet $X$ is understood.

Trimness and the preservation of the start state imply that there is at most one morphism between any two automata $(I_1, Q_{1},X, \square_{1})$ and 
$(I_2, Q_{2},X, \square_{2})$. Indeed, by trimness every state $q \in Q_{1}$ is of the form $q = I_{1}w$ for some $w \in X^{*}$ and thus if there is a morphism 
$f \colon Q_{1} \rightarrow Q_{2}$, we must have that $qf=(I_{1}w)f=(I_{1}f)w=I_{2}w$. Thus there is a morphism if and only if whenever $I_{1}w=I_{1}v$ 
for $v,w \in X^{+}$, then $I_{2}w =I_{2}v$ and this morphism is unique.

Let $\mathsf{A}_{1},\mathsf{A}_{2}$ be automata in $\mathcal{A}(X)$. It follows that if we define $\mathsf{A}_{1} \leqslant \mathsf{A}_{2}$ to 
mean that there is a surjective morphism from $\mathsf{A}_{2}$ to $\mathsf{A}_{1}$, then we can consider $\leqslant$ to be a partial order on $\mathcal{A}(X)$. 
We show now that this poset is a lattice. 

Let $\mathsf{A}_{1} = (I_1, Q_{1},X, \square_{1})$ and $\mathsf{A}_{2}= (I_2, Q_{2},X, \square_{2})$ be automata in $\mathcal{A}(X)$. We define 
$\mathsf{A}_{1} \vee \mathsf{A}_{2}$ to be the automaton with states $Q = \{(I_{1}w,I_{2}w) \mid w \in X^{*}\}$, initial state $(I_{1},I_{2})$ and action 
$(q_{1},q_{2})x=(q_{1}x,q_{2}x)$ for $(q_{1},q_{2}) \in Q$ and $x \in X$. The projections $\pi_{i}\colon Q \twoheadrightarrow Q_{i}$ for $i=1,2$ that send 
$(q_{1},q_{2})$ to $q_{i}$ are surjective morphisms. It is routine to check that $\mathsf{A}_{1} \vee \mathsf{A}_{2}$ is the join of $\mathsf{A}_{1}$ and 
$\mathsf{A}_{2}$ in the poset $(\mathcal{A}(X),\leqslant)$. Therefore $(\mathcal{A}(X),\leqslant)$ is a join semilattice. Since the $X$ automaton $1_{X}$ with 
one state is the minimal element of $(\mathcal{A}(X),\leqslant)$ and each principal down ideal is finite (by finiteness of states), we can define the meet 
of two automata to be the join of all elements less than or equal to each of them as usual. Therefore $(\mathcal{A}(X),\leqslant)$ is a lattice and has finite 
principal downsets.

We take a closer look at the join in the case of Cayley graphs of $X$-generated semigroups. Let $\phi_{1}\colon X^{+} \twoheadrightarrow S_{1}$ and 
$\phi_{2}\colon X^{+} \twoheadrightarrow S_{2}$ be presentations of the semigroups $S_{1}$ and $S_{2}$. Let $\Gamma_{1}$ and $\Gamma_{2}$ be the 
corresponding right Cayley automata. We define $\phi\colon  X^{+} \rightarrow S_{1} \times S_{2}$ by $x\phi = (x\phi_{1},x\phi_{2})$ for all $x \in X$. We write 
$\operatorname{Im}(\phi) = S_{1} \vee S_{2}$. Clearly $\Gamma_{1} \vee \Gamma_{2}$ is the completion of the Cayley automaton of $S_{1} \vee S_{2}$.

Note that $S_{1} \vee S_{2}$ is a subdirect product of $S_{1}$ and $S_{2}$. We have the following diagram in $\mathcal{A}(X)$
\[
	 \Gamma_{1} \stackrel{\pi_{1}}{\twoheadleftarrow \!\!\!-\!\!-} \Gamma_{1} \vee \Gamma_{2} \stackrel{\pi_{2}}{-\!\!-\!\!\!\twoheadrightarrow} \Gamma_{2}
\]
and corresponding morphisms in the category of $X$-generated finite semigroups
\[
	S_{1} \stackrel{\pi_{1}}{\twoheadleftarrow \!\!\!-\!\!-} S_{1} \vee S_{2} \stackrel{\pi_{2}}{-\!\!-\!\!\!\! \twoheadrightarrow} S_{2},
\]
where we use the same name for the corresponding projection morphisms. In the language of semigroup theory~\cite{qtheory}, we have the associated 
onto relational morphism $r= \pi_{1}^{-1}\pi_{2}\colon S_{1} \twoheadrightarrow S_{2}$. The following fact is straightforward, but important. We call this 
{\em string computation} of the relational morphism $r$.

\brm
\label{stringcomp1}
For all $s_{1} \in S_{1}$ we have $s_{1}r=\{I_{2}w \mid w \in X^{+}, I_{1} w=s_{1}\}$.
\erm

The notion of relational morphism can be generalized to the whole category $\mathcal{A}(X)$. Thus if $\mathsf{A}_{1} = (I_1, Q_{1},X, \square_{1})$ and
$\mathsf{A}_{2} = (I_2, Q_{2},X, \square_{2})$ are complete yautomata we have the following diagram:
\[
	Q_{1} \stackrel{\pi_{1}}{\twoheadleftarrow \!\!\!-\!\!-} Q_{1} \vee Q_{2} \stackrel{\pi_{2}}{-\!\!-\!\!\!\! \twoheadrightarrow} Q_{2}.
\]
Here $\pi_{1}, \pi_{2}$ are the projections and $r_{\mathsf{states}}= y\pi_{1}^{-1}\pi_{2}\colon Q_{1} \twoheadrightarrow Q_{2}$ is what we term a relational morphism 
between $\mathsf{A}_{1}$ and $\mathsf{A}_2$. Let $S_{i} $ for $i = 1,2$ be the transition semigroups of $\mathsf{A}_{i}$. Then one has a companion 
relational morphism $r_{\mathsf{sgp}} \colon S_{1} \rightarrow S_{2}$. See \cite{qtheory}. We then obtain a corresponding diagram of morphisms of 
transformation semigroups:
\[
	(Q_{1},S_{1}) \stackrel{\pi_{1}}{\twoheadleftarrow \!\!\!-\!\!-} (Q_{1} \vee Q_{2},S_{1} \vee S_{2}) 
	\stackrel{\pi_{2}}{-\!\!-\!\!\!\! \twoheadrightarrow} (Q_{2},S_{2}).
\]
In this context, Fact~\ref{stringcomp1} becomes the following.

\brm
\label{stringcomp2}
Let $r_{\mathsf{states}}\colon Q_{1} \rightarrow Q_{2}$ be the relational morphism defined above. Then 
\[
	(q_{1})r_{\mathsf{states}}=\{I_{2}w \mid w \in X^{*}, I_{1}w = q_{1}\}.
\]
\erm

%%%%%%%%%%%%%%%%%%%%%%%%%%%%%%%%%%%%%%%%%%%%%%%%%%%%%%%%%%%%%%%%%%%%%

\subsection{Preflows, partial flows and their relation to flows}
\label{partflows}

Let $\mathsf{A}=(I,Q,X,\square)$ be an object of $\mathcal{A}$. Recall that we have a fixed $X$-generated ${\rm GM}$ semigroup $(S,X)$, with 0-minimal 
ideal $M^{0}(A,G,B,C)$ and that $Rh_{B}(G)$ is the Rhodes lattice with $B$ and group $G$. We define a {\em preflow} on $\mathsf{A}$ to be any 
function $h\colon Q \rightarrow Rh_{B}(G)$. A preflow is called a {\em partial flow} on $\mathsf{A}$ if for all 
$p,q \in Q, x \in X$, if $px \neq\square$, then $(pf) \stackrel{\cdot x}{\rightarrow}(pf)x \neq\square$ satisfies the flow condition. We emphasize that partial flows are not defined on $\square$. Note that the restriction of a flow to $Q$ is a partial 
flow and a partial flow with $Q$ an invariant set in $\mathsf{A}$ is a flow.

As above, we have a fixed $X$-generated ${\rm GM}$ semigroup $(S,X)$. Let $\mathsf{A}=(I,Q,X,\square)$ be an  automaton with aperiodic 
transition semigroup $S=S(\mathsf{A})$. Let $m$ be the least integer such that $s^{m}=s^{m+1}$ for all $s \in S$. Recall that if $w \in X^+$, then 
$\mathsf{str}_m(w)$ is the sub-automaton of the Cayley graph of $B_{X}(w)$ whose states are all the elements that are greater or equal to $w\theta$ in the $\mathcal{R}$-order with initial state $I$ and terminal state $w\theta$ where $\theta$ is the natural morphism from $X^+$ to $B_{|X|}(m)$. This automaton accepts the language of words that have the same value as $w$ in the free Burnside semigroup $B_{|X|}(m)$. Its completion, which for ease we also denote by $\mathsf{str}_m(w)$ is a member of the category $\mathcal{A}$.

The following fact follows easily from the definitions.

\brm
 \label{flowlift}
Let $\mathsf{A}$ and $\mathsf{B}$ be automata over the same alphabet $X$. If 
$h\colon\mathsf{B} \rightarrow \mathsf{A}$ is an automaton morphism and $f\colon\mathsf{A} \rightarrow Rh_{B}(G)$ is a (partial) flow, then 
$hf\colon \mathsf{B} \rightarrow Rh_{B}(G)$ is a (partial) flow. We say that $h$ lifts the (partial) flow $f$.
\erm

\bl
\label{lemma.partial}
Let $\mathsf{A}=(I,Q,X,\cdot)$ be an aperiodic automaton with exponent $m$. If there is a flow $f \colon Q \rightarrow Rh_{B}(G)$, then 
the automaton $\mathsf{str}_{m_0}(w)$ has a partial flow for all $w \in X^{+}$ and for all $m_0 \geqslant m$.
\el

\begin{proof}
Recall that the non-sink states of $\mathsf{str}_{m_0}(w)$ are the elements of $B_{|X|}(m_0)$ that are $\mathcal{R}$ greater or equal to $w\theta$, where 
$\theta\colon X^{+} \rightarrow B_{|X|}(m_0)$ is the canonical morphism. It follows that if  
$u,v \in X^{+}$ are such that $u\theta = v\theta \geqslant_\mathcal{R} w\theta$, then $I'u=I'v$, where $I'$ is the start state of $\mathsf{str}_{m_0}(w)$.
Since there is a natural morphism $\alpha\colon B_{|X|}(m_0) \rightarrow S(\mathsf{A})$, it follows that the map sending $u$ to $Iu$ is a well-defined 
morphism in $\mathcal{C}$ and therefore by Fact~\ref{flowlift}, $\alpha f$ is a partial flow on $\mathsf{str}_{m_0}(w)$. 
\end{proof}

Let $s \in B_{n}(m)$ and $\theta \in Rh_{B}(G)$. Consider the preflow $f_{\theta}$ on 
$\mathsf{str}_{m}(s)$ that sends $I$ to $\theta$ and all other states 
to the empty $\mathsf{SPC}$, that is, to the bottom of the lattice $Rh_{B}(G)$. Then there is at most one minimal partial flow $f$ on $\mathsf{str}_{m}(s)$ 
such that $If = \theta$. This is because the collection of partial flows is closed under meet. We remark that there may be no partial flow extending this preflow, if 
the preflow forces us to arrive at the contradiction $\Longrightarrow\Longleftarrow$. Recall the definition of the set of $\mathsf{States}$ from 
Section~\ref{Ets}. We will show that the following Theorem is equivalent to Theorem~\ref{MainTh}. 

\bt
\label{equivMainTh}
Let $(S,X)$ be an $X$-generated ${\rm GM}$ semigroup with $\mathsf{RLM}(S)c = 1$. 
Then there is a computable $m>0$ such that $Sc = 1$ if and only if for each $\theta \in \mathsf{States}$, the preflow $f_{\theta}$ extends to a partial flow in $\mathsf{str}_{m}(s)$ 
for all $s \in B_{|X|}(m)$. Furthermore it is decidable for this value of $m$ if each $\mathsf{str}_{m}(s)$ has such a partial flow.
\et

\subsubsection{Bottom up definition of flow on a fixed automata}
\label{section.bottom up flow}

We now define an iterative procedure to compute the partial flow associated to an initial preflow. The  {\em initial preflow} $\ell_1$ is defined by
\[
	(I) \ell_1 = b_0 , (q)\ell_1 = \emptyset, q \neq I. 
\]
Here $b_{0} \in B$ is a fixed element of our ${\rm GM}$ semigroup $(S,X)$ with 0-minimal ideal $M^{0}(G,A,B,C)$. We identify $b \in B$ with the 
$\mathsf{SPC}$ $b/<1>$. Also $\emptyset$ represents the empty $\mathsf{SPC}$ which is the bottom of $Rh_{B}(G)$. 

We now iterate the initial preflow under two moves:
\begin{enumerate}
\item \textbf{Move forward:} At $q$, change $(q\cdot x)\ell$ to $(q)\ell \cdot x \vee (q\cdot x)\ell$.
\item \textbf{Move flow back under $\cdot x$:} At $q$, if $(q)\ell$ has blocks $B_1,\ldots, B_k$, all $B_j \cdot x \neq \emptyset$ and all
$B_1 \cdot x, \ldots, B_k \cdot x$ form a block of $(q \cdot x) \ell$, then consider the union $B_1 \cup \cdots \cup B_k$ and put on the 
unique cross section so that $\cdot x$ satisfies the flow condition. This will be the new $(q \cdot x) \ell$.
\end{enumerate}

We remark that in part (2) above, the unique cross section comes from the Tie-Your-Shoe Lemma~\cite[Chapter 4]{qtheory}. We also note that 
the automaton $\mathsf{A}$ has a unique pointwise minimal flow $f\colon I \rightarrow Rh_{B}(G)$ with $I(f)$ = $b_0$ or has no flow, since the 
collection of flows with this initial condition is closed under meet.

If $\ell_j$ is the preflow at the $j$-th iteration, then clearly $\ell_j \leqslant \ell_{j+1}$ pointwise. So either $\ell_{j+1}$ has the contradiction $\widehat{1}$
as a value or $\ell_j = \ell_{j+1}$ is a flow (which is the unique minimal flow with base point $b_0$). Indeed, if $\widehat{1}$ is never reached, then
$\ell_j=\ell_{j+1}$ is a flow since clearly $(q)\ell \cdot x \leqslant (q \cdot x) \ell$ by the forward move which is a flow by the flow back under $\cdot x$.

%%%%%%%%%%%%%%%%%%%%%%%%%%%%%%%%%%%%%%%%%%%%%%%%%%%%%%%%%%%%%%%%%

\subsection{Profinite methods}

We begin with an updated version of~\cite[Proposition 7.2]{stablepairs}. We first quote that proposition.

\bp\label{Trans7.2} \cite[Proposition 7.2]{stablepairs}
Let $X$ be a finite set and let $\alpha, \beta, \gamma$ be elements of $\widehat{F_{Ap}}(X)$ the free pro-finite aperiodic monoid on $X$. Then
$\alpha\beta\gamma = \alpha\beta$ if and only if one of the following three situations occur:
\begin{enumerate}

\item \textbf{Profinite separation:} $\beta\gamma = \beta$.
\item \textbf{Roll back I:} There exists $\tau \in \widehat{F_{Ap}}(X)$ such that $\alpha\beta\tau = \alpha$ and $\gamma = \tau\beta$.
\item \textbf{Roll back II:} There exist $\sigma, \tau \in \widehat{F_{Ap}}(X)$ and $i \geqslant 1$ such that  $\alpha\tau\sigma=\alpha, \beta = (\tau\sigma)^{i}\tau$,
and $\gamma = \sigma\tau$.
\end{enumerate}
\ep
\noindent

A change of variables shows that Condition 2 of Proposition \ref{Trans7.2} is the case $i=0$ of Condition~3. Indeed in Condition~3, if we take $i=0$ we have
$$\alpha\tau\sigma = \alpha, \beta = \tau, \gamma = \sigma\tau.$$
So we can eliminate $\tau$ by substituting $\beta$ and we have:
$$\alpha\beta\sigma = \alpha,  \gamma = \sigma\beta$$
and this is precisely Condition 2 of Proposition \ref{Trans7.2} with $\sigma$ playing the role of $\tau$ in that statement.

We therefore have a  version of~\cite[Proposition 7.2]{Trans} in which there are only two conditions.

\bp\label{New7.2}
Let $X$ be a finite set and let $\alpha, \beta, \gamma$ be elements of $\widehat{F_{Ap}}(X)$ the free pro-finite aperiodic monoid on $X$. Then
$\alpha\beta\gamma = \alpha\beta$ if and only if one of the following two situations occur:
\begin{enumerate}
\item \textbf{Profinite separation:} $\beta\gamma = \beta$.
\item  \textbf{Roll back:} There exist $\sigma, \tau \in \widehat{F_{Ap}}(X)$ and $i \geqslant 0$ such that  $\alpha\tau\sigma=\alpha, \beta = (\tau\sigma)^{i}\tau$,
and $\gamma = \sigma\tau$.
\end{enumerate}
\ep

Note that~\cite[Proposition 7.2]{Trans} does not say this specifically, but describes the right Stabilizer $\operatorname{Stab}(\alpha\beta)$ in $\widehat{F_{Ap}}(X)$. Part 1 of Proposition~\ref{New7.2} is the obvious 
$\operatorname{Stab}(\beta) \subseteq \operatorname{Stab}(\alpha\beta)$. Part 2 of Proposition~\ref{New7.2} relates $\operatorname{Stab}(\alpha\beta)$ to $\operatorname{Stab}(\alpha)$. Let us be more specific about this.

Recall that in a semigroup $S$ elements $s,t$ are conjugate if there are $x,y \in S$ such that $s=xy,t=yx$. If $G$ is a group, then two elements are conjugate in this sense if and only if they are conjugate in the usual sense of group theory. Two subsemigroups $T,U$ of $S$ are conjugate if there are elements $x,y \in S$ such that $xy \in T, yx \in U$ and $xUy \subseteq T, yTx \subseteq U$. A straightforward exercise shows that if $H,K$ are subgroups of a group $G$, then $H,K$ are conjugate in the above sense if and only if they are conjugate in the usual sense of group theory. 

In Statement 2 of Proposition~\ref{New7.2} note:

\begin{enumerate}
\item{$\gamma$, an element of $\operatorname{Stab}(\alpha\beta)$, is conjugate (in the above sense) to an element of
$\operatorname{Stab}(\alpha)$.}
\item{An easy induction on $i$ shows that $\alpha\tau = \alpha\beta$ and thus $\alpha\beta\sigma=\alpha$. That is, we have 
$\alpha\beta\mathcal{R}\alpha$. It also easily follows that $\operatorname{Stab}(\alpha\beta)$ is conjugate to 
$\operatorname{Stab}(\alpha)$ in this case.}
\end{enumerate}

So we can rewrite Proposition~\ref{New7.2} as follows giving an updated statement of~\cite[Proposition 7.2]{Trans}.

\bp \label{Final7.2}
Let $X$ be a finite set and let $\alpha, \beta, \gamma$ be elements of $\widehat{F_{Ap}}(X)$ the free pro-finite aperiodic monoid on $X$. Then
$\alpha\beta\gamma = \alpha\beta$ if and only if one of the following two situations occur:
\begin{enumerate}
\item \textbf{Profinite separation:} $\beta\gamma = \beta$.
\item  \textbf{Roll back:} $\alpha\beta\mathcal{R}\alpha$ and there are conjugate elements $\tau\sigma \in \operatorname{Stab}(\alpha), \sigma\tau \in 
\operatorname{Stab}(\alpha\beta)$ and an $i \geqslant 0$ such that $\beta = (\tau\sigma)^{i}\tau$,
and $\gamma = \sigma\tau$. In particular, $\operatorname{Stab}(\alpha\beta)$ is conjugate to $\operatorname{Stab}(\alpha)$.
\end{enumerate}
\ep

\bc
Let $X$ be a finite set and let $\alpha, \beta$ be elements of $\widehat{F_{Ap}}(X)$ the free pro-finite aperiodic monoid on $X$. Then one of the following two possibilities hold.
\begin{enumerate}

\item{$\operatorname{Stab}(\alpha\beta) = \operatorname{Stab}(\beta)$.}

\item{ $\alpha\beta\mathcal{R}\alpha$ and $\opn{Stab}(\alpha\beta)$ is conjugate to $\opn{Stab}(\alpha)$.}
\end{enumerate}
\ec

\brm

If $\alpha\beta$ is not $\mathcal{R}$-related to $\alpha$ then $\operatorname{Stab}(\alpha\beta) = \operatorname{Stab}(\beta)$ since $\operatorname{Stab}(\alpha\beta)$ can not be conjugate to $\opn{Stab}(\alpha)$ and thus Statement 2 of Proposition \ref{Final7.2} does not hold. The converse is not true. For example, if $\alpha=\beta$ is an idempotent, then obviously $\operatorname{Stab}(\alpha\beta) = \operatorname{Stab}(\beta)$ but $\alpha\beta = \alpha$.

\erm

In general, if $S$ is a semigroup and $s \in S$, then $\operatorname{Stab}(s)$ is the Stabilizer of $s$ in the right-Sch\"utzenberger representation of $S$ relative to the $\mathcal{R}$-class of $s$ and the conjugacy class of $\operatorname{Stab}(s)$ is in 1-1 correspondence with the $\mathcal{R}$-class of $s$. More general than this, stabilizers of points in transitive representations of $S$ are in 1-1 correspondence with the points of the representation. This is the generalization to semigroup theory that if $N$ is a subgroup of a group $G$, then the conjugacy class of $N$ in $G$ is in 1-1 correspondence with the cosets of $N$ in $G$, i.e. the points of the transitive permutation representation of $G$ on $G/N$ and each conjugate is the stabilizer of a unique coset. 

It is known that a stabilizer $\operatorname{Stab}(s)$ for $s \in \widehat{F_{Ap}}(X)$ is an internal $\mathcal{L}$-class meaning that any pair of elements 
has the property that one of them is a left multiple of the other by some element of $\operatorname{Stab}(s)$. Since it is known that stabilizers are 
closed subsemigroups in $\widehat{F_{Ap}}(X)$ this implies that $\operatorname{Stab}(s)$ has a minimal ideal $I(\operatorname{Stab}(s))$. The 
internal $\mathcal{L}$ property implies that $I(\operatorname{Stab}(s))$ is a left-zero semigroup. It follows that if $c \in I(\operatorname{Stab}(s))$ and 
$c' \in \operatorname{Stab}(s)$ then $cc'=c$.

Using the notation in Proposition \ref{Final7.2}, assume now that $c \in I(\operatorname{Stab}(\alpha\beta))$ is actually in 
$\operatorname{Stab}(\beta)$, that is, it satisfies $\beta c = \beta$. Then for all $c' \in \operatorname{Stab}(\alpha\beta)$ we have
$\beta c'= \beta cc'=\beta c=\beta$. Therefore, $\operatorname{Stab}(\alpha\beta)=\operatorname{Stab}(\beta)$ if and only if 
some (equivalently all) element(s) $c \in I(\operatorname{Stab}(\alpha\beta))$ is in $\operatorname{Stab}(\beta)$. The impact of this on what we are 
doing is that either we are already in Case 1 what we call {\em profinite separation} if $\operatorname{Stab}(\alpha\beta)=\operatorname{Stab}(\beta)$ 
or we can ``Rollback'' (and now there is only 1 Rollback) to a conjugate of $\operatorname{Stab}(\alpha)$ otherwise. So up to changing where we 
enter an $\mathcal{R}$-class when reading a word in a Cayley graph we can assume we are in Case 1.

We now apply the methods of \cite{stablepairs} to Proposition~\ref{Final7.2} to obtain a new version of~\cite[Corollary 7.3]{stablepairs} 
that we will use extensively in this paper. We first recall some details from \cite{Trans}.

Let $X$ be a finite set, $\widehat{X^{+}}$ the free profinite semigroup generated by $X$ and as above, $\widehat{F_{\opn{Ap}}}(X)$ the free pro-aperiodic 
semigroup generated by $X$. If $S$ is an $X$-generated semigroup then there are canonical morphisms $\alpha \colon \widehat{X^{+}}\rightarrow S$ 
and $\beta \colon \widehat{X^{+}}\rightarrow \widehat{F_{\opn{Ap}}}(X)$. The canonical relational morphism $\rho:S \rightarrow \widehat{F_{\opn{Ap}}}(X)$ 
is given by $\rho=\alpha^{-1}\beta$. For more details on the canonical relational morphism, see~\cite[Section 2]{stablepairs}. In that paper, there are canonical 
relational morphisms relative to an arbitrary pseudvariety, but we only consider the pseudovariety of aperiodic semigroups in this paper.

The following result summarizes~\cite[Theorem 2.6 and Lemma 4.4]{stablepairs} which are used crucially in that paper as it will be in the current paper.

\bl\label{discontmorphism}
\begin{enumerate}
\item{Let $S$ be an $X$ generated semigroup and let $\rho \colon S \rightarrow \widehat{F_{\opn{Ap}}}(X)$ be the canonical relational morphism. 
Then for all $\alpha \in \widehat{F_{\opn{Ap}}}(X)$, we have that $\alpha\rho^{-1} \in \mathsf{Pl}_{\mathsf{Ap}}(S)$. Every maximal pointlike set is of this form. 
Furthermore, the map $\rho^{-1} \colon \widehat{F_{\opn{Ap}}}(X)\rightarrow \mathsf{Pl}_{\mathsf{Ap}}(S)$ is a morphism of semigroups. It is not 
necessarily a continuous map between profinite semigroups.}
\item{Let $\alpha \in \widehat{F_{\opn{Ap}}}(X)$. Then $(\alpha,\operatorname{Stab}(\alpha))\rho^{-1}$ is a stable pair. Furthermore every maximal stable pair 
is of this form.}
\end{enumerate}
\el

We note that stable pairs are described in~\cite[Section 6]{stablepairs}, so we can compute all the maximal stable pairs.  
We need the following extension of~\cite[Theorem 6.4]{stablepairs}.

\bt \label{minidstab}
Suppose that $M$ is a finite monoid. Then the maximal 
stable pairs of $M$ are the maximal pairs (with respect to inclusion) $(Y,N)$ such that $Y \in \mathsf{Pl}_{\mathsf{Ap}}(M)$ and
there exists a submonoid $W$ of $\mathsf{Pl}_{\mathsf{Ap}}(M)$ with $W$ an internal $\mathcal{L}$--chain and:
\begin{enumerate}

\item{$\bigcup_{w \in W} W = N$ }

\item{$W$ is a subsemigroup of $Stab(Y)$ the pointwise stabilizer of $Y$ considered as an element of $\mathsf{Pl}_{\mathsf{Ap}}(M)$. }

\item{The minimal ideal of $W$ is a left-zero semigroup.}
\end{enumerate}
\et

\begin{proof}
The first two statements are precisely Theorem 6.4 of \cite{stablepairs}. By Lemma \ref{discontmorphism} there is an element $\gamma$ in the free profinite monoid $\widehat{F_{\opn{Ap}}}(X)$ generated by $X$ where $X$ is a fixed generating set of $M$ such that $Y=\gamma\rho_{\opn{Ap}}^{-1}$ and $N=\opn{Stab}(\gamma)\rho_{\opn{Ap}}^{-1}$ where $\rho_{\opn{Ap}}$ is the canonical relational morphism $M \rightarrow \widehat{F_{\opn{Ap}}}(X)$. Furthermore by Lemma \ref{discontmorphism} $f_{\opn{Ap}}=\rho_{\opn{Ap}}^{-1}$ is a not necessarily continuous morphism 
$f_{\opn{Ap}} \colon \widehat{F_{\opn{Ap}}}(X) \rightarrow \opn{Pl_{Ap}}(M)$. The proof of Lemma 6.4 of \cite{stablepairs} shows that $W=\widehat{F_{\opn{Ap}}}(X)f_{\opn{Ap}}$ is the $W$ in the statement of the Theorem.

Since $\opn{Stab}(\gamma)$ is a closed submonoid of $\widehat{F_{\opn{Ap}}}(X)$ and an internal $\mathcal{L}$-chain it follows that the minimal ideal $I$ of $\opn{Stab}(\gamma)$ is a left-zero semigroup. Thus for all $m \in I, n \in \opn{Stab}(\gamma)$, we have $mn=m$. It follows that every element of $I'=If_{\opn{Ap}}$ is stabilized on the right by every element of $W=\opn{Stab}(\gamma)f_{\opn{Ap}}$. It follows that $I'$ is the minimal ideal of $W$ and is a left-zero semigroup.
\end{proof}

\bc \label{corollary.stable}
Let $M$ be a finite monoid, $A,B \in \mathsf{Pl}_{\opn{Ap}}(M)$ with $N$ and $W$ as in Theorem \ref{minidstab}. Then the maximal stable pairs  
$(AB,N)$ are such that either $\operatorname{Stab}(AB) = \operatorname{Stab}(B)$ or are such that $A\mathcal{R}AB$ and 
$W$ is conjugate to a submonoid of $\operatorname{Stab}(A)$. 
\ec

We will use the notion of an inevitable graph from Section 2 of \cite{stablepairs}. We review the basic details from that paper.  

A directed graph $\Gamma$ consists of
a vertex set $V(\Gamma)$, an edge set $E(\Gamma)$ and functions $\iota, \tau: E(\Gamma) \rightarrow V(\Gamma)$
selecting the initial and terminal vertices $e$ of an edge, respectively. We
consider only finite graphs. A pointed graph is a graph together with a distinguished vertex $I$. 
An important example is the (right) Cayley graph of a semigroup $S$ with generating set $X$. This is the graph with vertex set $S^{I}=S \cup \{I\}$, where $I$ is an externally added identity element to $S$ and edge set $S^{I} \times X$ with $(s.x)\iota = s, (s,x)\tau=sx$ for all $s \in S^{I}, x \in X$.

A labeling of a graph $\Gamma$ over a transformation semigroup  $(X,S)$ is a pair of functions $\lambda =(\lambda_{V},\lambda_{E})$ $\lambda_{V}: V(\Gamma)  \rightarrow P(X), \lambda_{
E}:E(\Gamma)\rightarrow P(S)$, where $P(Q)$ is the power set of the set $Q$. If the image of $\lambda_{V}$ is contained in the singletons of $P(X)$ and the image of $\lambda_{E}$ is contained in the singletons of $S$, we call $\lambda$ a singleton labeling. A singleton labeling $\lambda$ is said to commute if,
for each edge $e$, we have $(e\iota\lambda_{V})(e\lambda_{E})=e\tau\lambda_{V}$. If $\phi\colon (X,S)\rightarrow (Y,T)$ is a relational morphism of transformation semigroups (\cite[Definition 4.14.8]{qtheory}), $\lambda$ is
a labeling of $\Gamma$ over $(X,S)$  and $\lambda'$ is a singleton labeling of $\Gamma$ over $(Y,T)$, then $\lambda$
is said to be $\phi$-related to $\lambda'$ if $v\lambda_{V} \subseteq (v\lambda'_{V})\phi^{-1}$ for all $v \in V(\Gamma)$, $e\lambda_{E} \subseteq e\lambda'\phi^{-1}$ for all $e \in E(\Gamma)$. We always assume that a relational morphism $\phi\colon (X,S) \rightarrow (Y,T)$ is parametrized by the companion relation of $f_{\phi}$ of $\phi$. See \cite{qtheory}.
This agrees with the definitions in~\cite[Section 2]{stablepairs} in the case we work with relational morphisms between the Cayley graphs of $S$ and $T$, but we need this more general notion.
 
We now generalize the definition of inevitability as described in~\cite[Section 2]{stablepairs} from monoids to transformation semigroups. We give the definition for the pseudovariety of aperiodic monoids as that is the only case we use in this paper.

\bd

Let $(X,S)$ be a transformation semigroup. A labeling $\lambda$ of a graph $\Gamma$ over $(X,S)$ is inevitable if, for all
relational morphisms $\phi\colon (X,S) \rightarrow (Y,T)$ with $T$ aperiodic there is a singleton labeling
$\lambda'$ of $\Gamma$ over $(Y,T)$ that commutes and that is $\phi$-related to $\lambda$.
\ed 

In this paper we use semigroups and semigroup Cayley graphs so we modify the above definition as follows. If $\phi\colon S\rightarrow T$ is a relational morphism we extend it to a relational morphism, (which we also denote by $\phi$ for convenience) $\phi\colon S^{I}\rightarrow T^{I}$ by letting $I\phi=\{I\}$.

For example, using pointed graphs, a subset $X$ of a semigroup $S$ is pointlike if the graph with two vertices $v,w$ and one edge $e$ with $e\iota=v,e\tau=w$, with label $I$ at $v$ and label $X$ on $e$ and $w$ is inevitable. Let $Y$ be a subset of a semigroup $S$ and $T$ a subsemigroup. Let $(Y,T)$ is a stable pair if the graph with two vertices $\{v,w\}$, $v$ and an edge from $v$ to $w$ labeled by $Y$, $w$ labeled by $Y$ and $|T|$ loops labeled bijectively by elements of $T$ is inevitable.

%%%%%%%%%%%%%%%%%%%%%%%%%%%%%%%%%%%%%%%%%%%%%%%%
\subsection{The setup for the proof}
\label{section.setup}

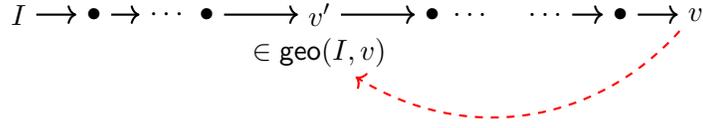
\begin{figure}[t]
\begin{center}
\begin{tikzpicture}[auto]
\node (A) at (0, 0) {$I$};
\node (B) at (1,0) {$\bullet$};
\node (CC) at (2,0) {$\cdots$};
\node (CCC) at (2.5,0) {$\bullet$};
\node (C) at (4,0) {$v'$};
\node(C0) at (4,-0.5)  {$\in \mathsf{geo}(I,v)$};
\node (D) at (5.5,0) {$\bullet$};
\node (DDD) at (6,0) {$\ldots$};
\node(DD) at (7,0) {$\ldots$};
\node (E) at (8,0) {$\bullet$};
\node (F) at (9,0) {$v$};
\path (A) edge[->,thick] (B);
\path (B) edge[->,thick] (CC);
\path (CCC) edge[->,thick] (C);
\path (C) edge[->,thick] (D);
\path (DD) edge[->,thick] (E);
\path (E) edge[->,thick] (F);
\path (F) edge[->,thick, bend left = 40, red, dashed] (C0);
\end{tikzpicture}
\end{center}
\caption{Schematic sketch of part of the McCammond tree $\mathcal{T}$ with black tree edges and a global edge indicated in dashed red.
\label{figure.geo1}}
\end{figure}

Let $\mathsf{A}=(I,Q,X,\square)$ be a finite state complete automaton. Let $f \colon Q \rightarrow Rh_{B}(G)$ be a flow. Assume that we have applied the $\operatorname{gst}$ expansion to $\mathsf{A}$.
Then by~\cite{gst}, $\mathsf{A}$ has a unique directed spanning tree $\mathcal{T}$. While $\mathcal{T}$ contains all the vertices of $Q$, 
$\mathsf{A}$ may have additional edges $\mathbf{E}$. The edges in $\mathbf{E}$ 
\[
	v \stackrel{x}{\longrightarrow} v'
\]
have the property that  $v'$ lies on the geodesic from the root $q_0$ to $v$ of the tree $\mathcal{T}$. The {\em geodesic} from $I$ to $v$ is
denoted by $\mathsf{geo}(I,v)$. The edges in $\mathbf{E}$ are called {\em global edges}, see~\cite{gst}. See also Figure~\ref{figure.geo1}.

In this section, $(S,Y, \phi)$ denotes a GM semigroup $S$ together with a surjective morphism $\phi \colon Y^{+}\rightarrow S$. If the morphism $\phi$ 
is understood in context, we write $(S,Y)$.

Let $X=(\States, \eval(S))$ be the evaluation transformation semigroup. We define a surjective morphism 
$\alpha \colon X^{\gst} \twoheadrightarrow X$ by composing the covering map $X^{\operatorname{Mc} \circ \operatorname{KR}} 
=X^{\gst}\twoheadrightarrow X^{\operatorname{KR}}$ with the expansion map $X^{\operatorname{KR}} \rightarrow X$. We then 
construct a relational morphism
\[
	r' \colon X^{\gst} \twoheadrightarrow (P,T),
\]
where $(P,T)$ is an aperiodic transformation semigroup with generating set $\eval(S)$ and satisfying the following properties:
\begin{enumerate}
	\item [(a)] $r'$ computes pointlike subsets of $\eval(S)^{\opn{gst}}$. That is, for all 
	$s \in T$, $sr'^{-1} \in \mathsf{Pl}_{\mathsf{Ap}}(\eval(S)^{\gst})$.
	\item[(b)] $r'$ computes stable pairs of $\eval(S)^{\gst}$. That is, for all $s \in T$, 
	$(s,\mathsf{Stab}_{R}(s))r'^{-1}$ is a stable pair.
       \item[(c)] $r'$ computes the inevitability of the labeled pointed graph
       \begin{equation}
       \label{equation.graph}
\raisebox{-0.5cm}{
\begin{tikzpicture}[auto]
\node (Q) at (0, 0) {$I$};
\node (Q1) at (1.5,0) {$A$};
\node (Q2) at (3,0) {$AB$};

\draw[edge,thick] (Q) -- (Q1) node[midway, above] {$A$};
\draw[edge,thick,bend right] (Q1) -- (Q2) node[midway, above] {$B$};
\path (Q2) edge [thick,loop right] node {$C$} (Q2);
\end{tikzpicture}}
\end{equation}
where $A,B, \subseteq \eval(S)^{\opn{gst}} , C=\{c_{1},c_{2}, \ldots, c_{l}\}$ for some $l>0$. 
\end{enumerate}

Let $\alpha\colon (P,T)^{\gst} \twoheadrightarrow (P,T)$ be the covering map. We define
\begin{equation}
\label{equation.definition r}
	r\colon X^{\gst} \twoheadrightarrow (P,T)^{\gst}
\end{equation}
by $r=r'\alpha^{-1}$.
We write $(Q,A)$ for $(P,T)^{\gst}$. It is known that $A$ is an aperiodic semigroup~\cite{gst}. It follows from Proposition~\ref{proposition.properties} 
that $r$ satisfies the same properties (a) and (b) that $r'$ does. 

We know that it is decidable to effectively find such a $(Q,A)$ and relational morphism $r$. 
This makes the graph in~\eqref{equation.graph}
inevitable as in~\cite[Lemma 2.6]{stablepairs} and Corollary~\ref{corollary.stable}.

We have the relational morphism $r^{-1} \colon (Q,A,\eval(S)) \rightarrow (X^{\gst},\eval(S))$. This implies that 
for all $s \in A$, we have $sr^{-1}=\{I_{X^{\gst}}w \mid w \in \eval(S), I_{A}w=s\} \in \mathsf{Pl}_{\mathsf{Ap}}(\eval(S),
\eval(S))$, since we have assumed that $r$ computes the pointlike subsets of $\eval(S)$. 

\brm
\label{remark.aug}
\mbox{}
\begin{enumerate}
	\item For all $s \in A$, if $sr^{-1} \in \mathsf{Pl}_{\mathsf{Ap}}(\eval(S)) \cap S$ it is not necessarily in 
	$\mathsf{Pl}_{\mathsf{Ap}}(S)$.
	\item Since $sr^{-1} \in \mathsf{Pl}_{\mathsf{Ap}}(S)$ it follows that the join $\bigvee sr^{-1}$ over elements in $sr^{-1}$ is less than or equal to some (not necessarily unique) $f$ in $\eval(S)$ \cite[Lemma 4.28]{Trans}. We pick such an element $f$ and denote it by $\overline{\bigvee sr^{-1}}$.
    
    	The map that sends an element $Z$ of $\mathsf{Pl}_{\mathsf{Ap}}(\eval(S))$ to $\overline{\bigvee Z} \in \eval(S)$ is called the Augmentation map and will be studied in a later section
	(see Definition~\ref{definition.augmentation}).
	\item We will prove in the next sections that $(Q,A)$
	has a flow with no contradiction if and only if $Sc=1$ under the assumption that $\RLM(S)c \leqslant 1$. That is either $X^{\gst}$ has a flow with no contradiction and $Sc=1$ or no such flow exists 
	and $Sc=2$. This is the decidability criterion for complexity 1.
\end{enumerate}
\erm

%%%%%%%%%%%%%%%%%%%%%%%%%%%%%%%%%%%%%%%%%%%%%%%%%%
\subsection{The procedure}
\label{automaton1}

We now have the following relational morphism:
\[
\widehat{r}^{-1} \colon (Q,A,\Eval(S)) \twoheadrightarrow (X,\Eval(S)),
\]
that is the composition of $r^{-1}$ with the projection from $X^{\gst}$ to $X$. It is a relational morphism between pointed transformation semigroups, 
where $I$ is the initial point of $(Q,A,\eval(S))^{\gst}$ and we choose a fixed initial state $b_{0}/\langle 1\rangle$ for $(X,\eval(S))$ acting on $\States$. 

We have applied the gst expansion to $(P,T,\opn{Eval}(S))$ to obtain $(Q,A,\eval(S))$. We denote its unique directed spanning tree by $\mathcal{T}$.
We can pick a finite set of elements $\{t_{1}, \ldots , t_{k}\}$ of the free  semigroup $\eval(S)^{+}$ such that these elements map bijectively onto the tree 
$\mathcal{T}$. That is, $\mathcal{T}$ has $k$ leaves $q_{i}$ for $i=1, \ldots, k$ and each word $t_{i}$ for $i=1,\ldots, k$ reads a path in $\mathcal{T}$ 
from the initial state to $q_{i}$. 
Recall that if $t \in \Eval(S)^{+}$, then $\mathsf{str}_m(t)$ is the subautomaton of $(Q,A,\eval(S))$ consisting of all states that are on paths that start at $I_{A}$ and end at $I_{A}t$. Here $m$ is such that $t^m  =t^{m + 1}$. It follows that $(Q,A,\eval(S))$ is the join
$\bigvee_{i = 1}^{k}\mathsf{str}_m(t_i)$. 

Furthermore, each automaton $\mathsf{str}(t_{i})$ is a linear automaton \cite{gst}. That is, its strongly connected components are linearly ordered and its transition edges form a chain in the strongly connected component graph. We can thus factor each $t_{i}$ as $t_{i}=\tau_{i,1}u_{i,1}\tau_{i,2}u_{i,2}\ldots \tau_{i,j_{i}}u_{i,j_{i}}$ where the $\tau_{i,j}$ are precisely the transition edges of $t_{i}$ and some of the $u_{i,j}$ may be the empty word. Let the strongly connected components of $\mathsf{str}(t_{i})$ be $C_{i,0}=\{I\}, C_{i,1}, \ldots, C_{i,j_{i}}$. Then $\tau_{i,j}$ is the unique edge in $\mathsf{str}(t_{i})$ from $C_{i,j-1}$ to $C_{i,j}$. In the language of \cite{gst}, the prefix of $t_{i}$ ending in $\tau_{i,j}$ is the unique entry to $C_{i,j}$ and the prefix that ends in $u_{i,j}$ is the unique exit point from $C_{i,j-1}$ in $\mathsf{str}(t_{i})$.

We now define a procedure that defines a sequence of partial flows on $(Q,A,\eval(S))$. 
We will prove under the hypothesis of Theorem~\ref{equivMainTh} that these partial flows limit to a flow. 
Each partial flow $\phi_{i,j}$ will be defined on a subautomaton of $(Q,A,\eval(S))$. 
We begin by processing $t_1$. We use all the notation from the previous paragraph.

The domain of the initial partial flow is $\{I\}$. It sends $I$ to $b_{0}\space /\langle 1 \rangle$ where $(1,b_{0})$ is the designated start state in $G \times B$. By induction, the $j^{th}$ partial flow, $\phi_{1,j}$ for $1 \leqslant j < j_{1}$ has domain $D_{1,j}=C_{1,0} \cup \{\tau_{1,1}\} \cup \C_{1,1} \cup \cdots 
\cup\{\tau_{1,j}\} \cup C_{1,j}$. The domain of the next partial flow $\phi_{1,j+1}$ is
$D_{1,j} \cup \{\tau_{1,j+1}\} \cup C_{1,j+1}$. That is, we add the transition edge from $C_{1,j}$ to $C_{1,j+1}$ and the strongly connected component 
$C_{1,j+1}$ of the previous domain. We have the following Hypothesis.

\bh
\label{hyp1}
Let $q_{j}$ be the exit state from $C_{j}$ to $C_{j+1}$. Then the preflow $f_{(1,j+1)}$ on the subautomaton of $(Q,A)$ with domain 
$\{q_{j}\stackrel{\tau_{j+1}}{\rightarrow}q_{j}\tau_{j+1}\} \cup C_{j+1}$  that sends $q_{j}$ to $q_{j}\phi_{1,j}$ and all other states to the empty 
$\mathsf{SPC}$ extends to a minimal partial flow $\widetilde{f}_{1,j+1}$ with $q\widetilde{f}_{1,j+1} \geqslant q\phi_{1,j}$.
\eh

If Hypothesis \ref{hyp1} holds, then we obtain a partial flow $\phi_{(1,j+1)}$ with domain $D_{j+1}$ by letting it have the value of $\phi_{(1,j)}$ on 
$D_{j} -\{q_{j}\}$ and the value of $\widetilde{f}_{j+1}$ on 
$\{q_{j}\} \cup C_{1,j+1}$. By induction this procedure defines a partial flow on $\mathsf{str}(t_{1})$.

We continue the process with $\mathsf{str}(t_{2})$. If all the strongly connected components of $\mathsf{str}(t_{2})$ are contained in those of $\mathsf{str}(t_{1})$ then we proceed to $t_{3}$. Otherwise there exists a minimal transition edge (recall that the transition edges of $\mathsf{str}(t_{2})$ form a chain)  $\tau$ of 
$\mathsf{str}(t_{2})$ that transitions to a new strongly connected component. We now define a Hypothesis analogous to Hypothesis \ref{hyp1} that says that we can extend the current partial flow to one that adds the transition edge $\tau$ and the strongly connected component it enters to the domain. Throughout this discussion we use the fact that $(Q,A)$ has the unique simple path property, there is exactly 1 simple path from $I$ to any state $q \in Q$. In particular this implies that there is at most one transition edge between any two strongly connected components. See \cite{gst} for details.

\bh\label{mainhyp}

If the current partial flow processing $t_i$ does not contain all states in $Q$, then we can extend to a partial flow whose domain adds a transition edge and the strongly connected component it enters to its domain in analogy to the process in Hypothesis \ref{hyp1} and its conclusion.

\eh

It is clear that if we prove that Hypothesis \ref{mainhyp} holds then the process allows us to define a partial flow on all of $Q$, that is, to a flow. 

We have reduced proving that computing complexity 1 is decidable to proving that 
Hypothesis~\ref{mainhyp} is true for all 
$i$. In the next sections we do this by 
defining {\em witnesses}. A witness is a word $w \in \eval(S)^{+}$ together with the straightline automaton for the word problem $\mathsf{str}_m(w)$.
Assume that in Hypothesis \ref{mainhyp} we need to process the transition edge $\tau_{i,j}$. We will find a computable $m$ and a word $w_{i,j} \in \eval(S)^{+}$ so that the  preflow on $\mathsf{str}_m(w_{i,j})$ that sends $I$ to the value of the current partial flow's value on $\tau_{i,j}\iota$ and all other states to the empty $\mathsf{SPC}$ 
extends to a partial flow $f$ that witnesses the truth of Hypothesis \ref{mainhyp} for $\tau_{i,j}$.

%%%%%%%%%%%%%%%%%%%%%%%%%%%%%%%%%%%%%%%
\subsection{Profinite Witnesses} \label{profwitnesses}

We use all the notation from Section \ref{automaton1}. At each stage of the process leading to Hypothesis \ref{mainhyp}, we have a partial flow $\phi:D \rightarrow Rh_{G}(B)$ where $D$ is a union of strongly connected components of $(Q,A)$. If $D$ is all of $Q$, then $\phi$ is a flow. Otherwise, the process attempts to add a transition edge $\{q\stackrel{\tau}{\rightarrow}q\tau\}$ union the strongly connected component $C$ this edge enters to $D$ and extend the partial flow $\phi$ to this bigger domain. In this section, we first define witnesses to Hypothesis~\ref{mainhyp} in $\mathsf{Pl_{Ap}}(\eval(S))$ using profinite methods. In the next section we use 
McCammond's work~\cite{Mc91, Mc-nf} to turn them into witnesses in $\eval(S)$.

Using the relational morphism $\widehat{r}$ we let $\Theta = \tau\widehat{r}^{-1}$ and $W=\opn{Stab}(\tau)\widehat{r}^{-1}$. Then $(\Theta,W)$ is a 
stable pair by the choice of $\widehat{r}$. Since partial flows are closed under restriction of domain we can embed this stable pair into a maximal stable 
pair that we continue to denote by $(\Theta,W)$ and prove that the conditions of 
Hypothesis~\ref{mainhyp} hold. Therefore by Theorem \ref{minidstab}, $\Theta \in \mathsf{Pl_{Ap}}(\eval(S))$, $W \subseteq \opn{Stab}(\Theta)$ and is an
internal $\mathcal{L}$-chain whose minimal ideal is a left-zero semigroup. We call $(\Theta,W)$ a profinite witness for Hypothesis \ref{mainhyp}.

%%%%%%%%%%%%%%%%%%%%%%%%%%%%%%%%%%%%%%%%%%%%
\subsection{Witnesses in free $\omega$-aperiodic and free Burnside semigroups and the completion of the proof}
\label{section.witness}

In this section we let $n=|\eval(S)|$. We now use the results of~\cite{Mc91, Mc-nf} to turn the profinite witness of 
Section \ref{profwitnesses} into a witness of the form $\mathsf{str}_m(w)$ 
for some $w$ in the free aperiodic Burnside semigroup $B_{n}(m)$. We will determine $m$ in the proof and explain what we mean by ``witness" below.

We note that we have three interpretations of 
the $\omega+*$ operator on certain semigroups:

\begin{enumerate}
\item In the pointlike aperiodic semigroup $\mathsf{Pl_{Ap}}(S)$ of a semigroup $S$, if $Z \in \mathsf{Pl_{Ap}}(S)$, then 
$Z^{\omega+*} = Z^{\omega}(\cup_{i=1}^{\omega-1}Z^{i})$. Henckell's Theorem~\cite{Henckell} states that $\mathsf{Pl_{Ap}}(S)$ is the smallest 
subsemigroup of $P(S)$ containing the singleton sets, closed under subset and also the $\omega+*$ operator. See also Section~\ref{section.more tools}.
\item If $t \in \eval(S)$, then $t^{\omega+*}$ 
is defined as in~\cite{Trans}. Here we regard elements of $\eval(S)$ 
as operating on the Rhodes lattice $Rh_{B}(G)$.
\item In the free Burnside semigroup $B_{|X|}(m)$, $t^{\omega+*}=t^{\omega}$.
\end{enumerate}

\bd
\label{definition.augmentation}
We define the augmentation map
\[
	\mathsf{\mathsf{Aug}}\colon \mathsf{Pl_{Ap}}(\eval(S)) \rightarrow \eval(S)
\]
by $\mathsf{Aug}(Z) = \overline{\bigvee Z}$ (see Remark~\ref{remark.aug}).
\ed

 To prove Theorem \ref{MainTh} it suffices to prove that if there exists a flow then the procedure described in 
Section~\ref{automaton1} gives a method to obtain a flow on $(Q,A)$. We thus claim that if the procedure gives a flow then $Sc=1$. Otherwise the procedure leads to a contradiction, no flow exists and $Sc=2$.

We assume that a flow exists. Then we must prove that the procedure outlined in Section~\ref{automaton1} gives a flow.
By Lemma~\ref{lemma.partial} we can assume that there exists an $m_0$ such that $\mathsf{str}_{m}(t)$ has a partial flow for all $t$ and
all $m \geqslant m_0$. Note that the existence of $m_0$ is given but not necessarily its computability. However, the proof below gives computability.  

Now consider the procedure of Section~\ref{automaton1}. Assume that we are processing $t_i$ and that we wish to extend the currently defined partial flow $\phi$ by adding a transition edge $\{q\stackrel{\tau}{\rightarrow}q\tau\}$ and the strongly connected component $C$ that it enters to the domain of $\phi$. 

We begin with the profinite witness $(\Theta,W)$ that we constructed in the previous section.  By Definition \ref{definition.augmentation}, since  $\Theta \in \mathsf{Pl_{Ap}}(\eval(S))$ we have $\overline{\bigvee \Theta} \in \eval(S)$. 
We will need the concept of Well-Formed Formula (WFF) as defined in \cite[Definition 5.1]{Trans}. Let $\Delta$ be a WFF representing $\overline{\bigvee \Theta}$. It follows from the proof of~\cite[Proposition 4.28]{Trans} that we can use the same WFF $\Delta$ that represents $\Theta \in \mathsf{Pl}(\eval(S))$ 
for $\overline{\bigvee \Theta} \in \eval(S)$. Such a $\Delta$ for $\Theta$ 
exists by Henckell's Pointlike Theorem.

We wish to prove that $\overline{\bigvee \Theta} \in \eval(S)$ does not have the contradiction in its range as a function on $Rh_{B}(G)$. This is precisely 
what we need to show that we can extend our partial flow as per the procedure. By~\cite[Proposition 4.24]{Trans} we need to show that 
$\overline{\bigvee \Theta} \in \eval(S)$ has a value. See~\cite[Definition 4.23]{Trans}.  To do this, we will find a computable value of $m$ and an element 
$w \in B_{n}(m)$ that witnesses that $\overline{\bigvee \Theta} \in \eval(S)$ has a value in a sense to be made precise below.

We use the theory of free $\omega$-aperiodic and free aperiodic semigroups as in~\cite{Mc91,Mc-nf} to find a witness for Hypothesis~\ref{mainhyp} 
of the form $\mathsf{str}_m(t)$. We summarized the material we need in Section \ref{freeaperiodic}. The WFF $\Delta$ defines an element $t$ of the free $\omega$-semigroup $\Omega_{Ap}(n)$. Without loss of generality we can assume that $t$ is in normal form \cite{Mc-nf}. By the results of~\cite[Sections 12 and 13]{Mc-nf} (see also Theorem \ref{mainMcCam} and Section \ref{freeaperiodic}), there is a computable value $m_1$ such that the natural map from $\theta_{m}\colon \Omega_{Ap}(n)\rightarrow B_{n}(m)$ separates $t$ from all words of smaller length and for all $m\geqslant m_1$. In particular, $t$ is mapped by $\theta_{m}$ to the element with automaton  $\mathsf{str}_{m}(t\theta_{m})$ for $m\geqslant m_1$.

Let $m = \opn{max}\{m_{0},m_{1}\}$. As noted above, $\mathsf{str}_{m}(t)$ has a partial flow with no contradictions. Since $\overline{\bigvee \Theta} \in \eval(S)$ is defined by the same WFF as $t$ it follows that the contradiction is not in the range of $\overline{\bigvee \Theta} \in \eval(S)$ by~\cite[Proposition 4.24]{Trans}. Otherwise $\mathsf{str}_{m}(t)$ would not have a contradiction free partial flow. Furthermore, \cite[Proposition 4.24]{Trans} implies that $\overline{\bigvee \Theta} \in \eval(S)$ has a value $\nu$.

Let $q$ be the state we are processing in the procedure and $\phi$ the current partial flow. We can assume that $q\phi$ is not the empty $\mathsf{SPC}$. 
If $q\phi$ was empty we can extend the flow by sending $q\overline{\bigvee \Theta}$ and all states in the strongly connected component $C$ to be the empty $\mathsf{SPC}$ as well.

Now assume that in the procedure of Section~\ref{automaton1} both $q\phi$ is not empty and that the $(q\phi)\overline{\bigvee \Theta}$ is not empty.
Since by Theorem \ref{minidstab}, $W$ is a subsemigroup of $\opn{Stab}(\Theta)$, we have Profinite Separation $\Theta W=\Theta$ in $\mathsf{Pl_{Ap}(\eval(S))}$. This shows the importance of Profinite Separation.
Since $\Theta$ is pointlike, it is covered by a point in $\mathsf{str}_{m}(t)$ since this automaton has an aperiodic transition semigroup. Thus $\tau$ is readable in $\mathsf{str}_m(t)$ and since $\Theta W=\Theta$, $W$ lifts to a submonoid of the stabilizer of every point in the image of $\tau$ in $\mathsf{str}_{m}(t)$. Therefore $\tau$ is $\mathcal{J}$ above $\mathsf{red}(t)$, the representative of $t$ in $B_{n}(m)$ which by \cite{Mc91} is the unique shortest accepted by $\mathsf{str}_m(t)$ and the connected component $C$ of $q\tau$ lifts to a connected component of $\mathsf{str}_{m}(t)$. This is what we needed to prove and we are done. We now see that we can take $m_0$ to be the exponent of the graph $r$ in~\eqref{equation.definition r}. This completes the proof of Theorem \ref{MainTh}. 

\bibliography{stubib.bib}
\bibliographystyle{abbrv}

{(Stuart Margolis) Department of Mathematics, Bar Ilan University, Ramat Gan 52900, Israel

{\it Email address}: \; \texttt{margolis@math.biu.ac.il}}

\medskip

{(John Rhodes) Department of Mathematics, University of California, Berkeley, CA 94720, U.S.A.

{\it Email address}: \; \texttt{rhodes@math.berkeley.edu, blvdbastille@gmail.com}}

\medskip

{\large (Anne Schilling) Department of Mathematics, UC Davis, 
One Shields Ave., Davis, CA 95616-8633, U.S.A.

{\it Email address}: \; \texttt{aschilling@ucdavis.edu}}

% --------------------------------------------------------------

\end{document}